\documentclass[12pt]{article}

\usepackage[utf8]{inputenc}

\usepackage{tikz}
\usepackage{amsthm}
\usepackage{amsmath} 
\usepackage{amssymb} 
\usepackage{bm}
\usepackage{multirow}
\usepackage{booktabs}
\usepackage{hyperref}
\usepackage{caption}
\usepackage{subcaption}
\usepackage{algorithm}
\usepackage{algorithmic}

\graphicspath{ {./pic/}}

\renewcommand{\div}{\mathop{\rm div}\nolimits}

\makeatletter
\renewcommand{\@biblabel}[1]{#1.\hfill}
\makeatother

\numberwithin{equation}{section}

\theoremstyle{remark}
\newtheorem{remark}{\bf Remark}

\title{Meshfree GMsFEM-based exponential integration for multiscale 3D advection–diffusion problems}

\date{\today}

\author{Djulustan Nikiforov\thanks{North-Eastern Federal University, Yakutsk, Russia; \texttt{dju92@mail.ru}}  \and Leonardo A. Poveda\thanks{College of Science, Mathematics and Technology, Wenzhou-Kean University, Wenzhou, Zhejiang 325060, P. R. China; \texttt{lpovedac@wku.edu.cn}} \and  Dmitry Ammosov\thanks{Chemical and Petroleum Engineering Department, Khalifa University of Science and Technology, Abu Dhabi, 127788, UAE; \texttt{dmitrii.ammosov@ku.ac.ae}}  \and Yesy Sarmiento\thanks{Escuela de Matemáticas y Ciencias de la Computación, Universidad Nacional Autónoma de Honduras, Boulevard Suyapa, Tegucigalpa, M.D.C., Honduras; \texttt{yesy.sarmiento@unah.hn}}\and Juan Galvis\thanks{Departamento de Matemáticas, Universidad Nacional de Colombia, Carrera 45 No. 26-85, Edificio Uriel Gutiérrez, Bogotá D.C., Colombia; \texttt{jcgalvisa@unal.edu.co}}\and Mohammed Al Kobaisi\thanks{Chemical and Petroleum Engineering Department, Khalifa University of Science and Technology, Abu Dhabi, 127788, UAE; \texttt{mohammed.alkobaisi@ku.ac.ae}}
}

\begin{document}
\maketitle

\begin{abstract}
In this work, we extend the meshfree generalized multiscale exponential integration framework introduced in Nikiforov et al. (2025) to the simulation of three-dimensional advection--diffusion problems in heterogeneous and high-contrast media. The proposed approach combines meshfree generalized multiscale finite element methods (GMsFEM) for spatial discretization with exponential integration techniques for time advancement, enabling stable and efficient computations in the presence of stiffness induced by multiscale coefficients and transport effects. We introduce new constructions of multiscale basis functions that incorporate advection either at the snapshot level or within the local spectral problems, improving the approximation properties of the coarse space in advection-dominated regimes. The extension to three-dimensional settings poses additional computational and methodological challenges, including increased complexity in basis construction, higher-dimensional coarse representations, and stronger stiffness effects, which we address within the proposed framework. A series of numerical experiments in three-dimensional domains demonstrates the viability of the method, showing that it preserves accuracy while allowing for significantly larger time steps compared to standard time discretizations. The results highlight the robustness and efficiency of the proposed approach for large-scale multiscale simulations in complex heterogeneous media.
\end{abstract}

\section{Introduction}

Flow simulation in porous media plays a central role in many scientific and engineering applications, including groundwater management, oil reservoir simulation, subsurface carbon storage, and the design of porous materials \cite{narayanan2018flow,battiato2019theory, bastian2013simulation}. These problems are typically characterized by heterogeneous and high-contrast permeability fields, which introduce significant challenges for both spatial and temporal discretization. In particular, the permeability coefficient $\kappa(x)$ may exhibit strong variations across multiple scales, leading to severe computational demands when classical numerical methods are employed.

Modeling and simulation of such problems require addressing several key difficulties. First, resolving fine-scale variations of $\kappa$ often demands very fine spatial discretizations, making direct numerical simulation computationally prohibitive \cite{egh12, efendiev2011multiscale, abreu2019convergence, calo2016randomized}. Second, the presence of high-contrast coefficients introduces stiffness in time-dependent problems, restricting the stability region of classical time integration schemes such as Crank--Nicolson. These issues are further compounded in applications where the permeability is only available through indirect descriptions, such as black-box evaluations, parameterized models, or discrete data on complex meshes.

Generalized Multiscale Finite Element Methods (GMsFEM) have been developed to efficiently address the spatial complexity of multiscale high-contrast problems by constructing reduced-order approximation spaces based on local spectral problems \cite{chung2014generalized, gao2015generalized, fu2019generalized, chung2014adaptive}. These methods have been successfully applied to a wide range of problems, including flow, diffusion, and elasticity in heterogeneous media \cite{egh12, efendiev2011multiscale, abreu2019convergence, calo2016randomized, zambrano2021fast}. 

However, classical GMsFEM relies on the construction of coarse meshes, which can be computationally expensive or even impractical in complex geometries or when permeability data is defined on irregular grids. To address this limitation, meshfree variants of GMsFEM have been introduced. In particular, the meshfree GMsFEM framework developed in \cite{nikiforov2023meshfree} constructs coarse approximation spaces based solely on point clouds, avoiding the need for explicit coarse mesh generation. This approach has been extended to several classes of problems, including nonlinear and time-dependent settings \cite{djulustan2024meshfree, nikiforov2023meshfree,nikiforov2024meshfree1,nikiforov2025meshfree2}.

In our previous work \cite{nikiforov2025meshfree}, we introduced a meshfree generalized multiscale exponential integration method for parabolic problems. There, we combined meshfree GMsFEM for spatial discretization with exponential time integration to efficiently handle stiffness induced by high-contrast multiscale coefficients. The results demonstrated that this approach allows for stable and accurate simulations using significantly larger time steps compared to classical schemes, while maintaining reduced computational complexity through multiscale model reduction.

Despite these advances, many practical applications inherently require the simulation of fully three-dimensional problems. Extending meshfree GMsFEM combined with exponential integration to 3D settings is not a straightforward task. The transition from two-dimensional to three-dimensional problems introduces several additional challenges. First, the computational cost of constructing multiscale basis functions increases significantly due to the higher dimensionality of local spectral problems. Second, the number of degrees of freedom associated with both fine-scale and coarse-scale representations grows rapidly, affecting memory requirements and computational efficiency. Third, the geometric complexity of point cloud distributions in three dimensions complicates the construction of stable and accurate partition of unity functions. Finally, the stiffness of the resulting systems is typically more pronounced in 3D, making the design of efficient and stable time integration strategies even more critical.

The main objective of this work is to investigate the numerical viability of combining meshfree GMsFEM with exponential integration for three-dimensional advection--diffusion problems. In particular, we assess whether the advantages observed in two-dimensional settings---namely stability for large time steps and efficient multiscale approximation---can be preserved in 3D simulations. To this end, we extend the meshfree generalized multiscale exponential integration framework to three-dimensional problems and present a series of numerical experiments that evaluate its performance in terms of accuracy, stability, and computational efficiency.

Exponential integration methods have been shown to be effective for stiff problems arising in multiscale simulations \cite{contreras2023exponential, Expint2007,poveda2024second}. Alternative strategies, such as partially explicit time integration \cite{li2022partially, djulustan2024meshfree}, have also been explored. In this work, we focus on exponential integrators due to their ability to exploit the reduced-order structure provided by GMsFEM, enabling larger time steps without compromising stability.

The rest of the paper is organized as follows. Section \ref{sec:preliminaries} reviews the relevant background and existing methods, it also describes the exponential integration approach. Section \ref{sec:meshfreegmsfem} presents the meshfree GMsFEM framework including the 
GMsFEM version of the fine-grid time discretizations.  Section \ref{sec:numerical_results} presents numerical experiments in three-dimensional settings. Finally, Section 
\ref{sec:conclusions} concludes the paper.

\section{Problem formulation and fine-grid approximation}\label{sec:preliminaries}
In this section, we introduce the mathematical formulation of the time-dependent advection--diffusion problem considered in this work and describe its fine-grid discretization. The fine-scale model resolves the heterogeneous and high-contrast coefficients that characterize the multiscale nature of the problem and serves as the reference solution for evaluating the accuracy of the proposed method. We present the continuous formulation, followed by its spatial discretization on a sufficiently fine grid, and briefly discuss the associated computational challenges that motivate the use of reduced-order multiscale approaches.
\subsection{Problem formulation}


Let $\Omega \subset \mathbb{R}^3$ be the spatial computational domain and $I = [0, T]$ be the time domain. Then, we have the following advection-diffusion problem defined in $\Omega \times I$
\begin{equation}\label{eq:fine_model}
\begin{split}
\frac{\partial u}{\partial t} - \epsilon \div \left( \kappa(x) \nabla u \right) + \beta \cdot \nabla u = f(u) \quad &\text{in } \Omega \times I,\\
u = u_D \quad &\text{on } \partial \Omega \times I,\\
u(0, x) = u_0(x) \quad &\text{in } \Omega,
\end{split}
\end{equation} 
where 
$u$ is the solution field (e.g., concentration, temperature),
$\kappa(x)$ is the heterogeneous diffusion coefficient,
$\beta$ is the velocity field,
$\epsilon$ is a constant determining the ratio of advection and diffusion,
$f(u)$ is a source/sink term,
$u_D$ is the boundary value, 
and $u_0(x)$ is the initial value.
Note that we suppose that $\kappa(x)$ possesses a high contrast, i.e., $\text{max}_{x \in \Omega} / \text{min}_{x \in \Omega} \gg 1$.

\subsection{Fine-grid finite-element approximation}

Let us consider the fine-grid approximation of \eqref{eq:fine_model}, which we will use as the reference one in this work. First, we apply a finite element method with the standard piecewise linear basis functions for spatial approximation. Let us define the following trial and test function spaces 
\begin{equation*}
V := \{ v \in H^1(\Omega): v = u_D \text{ on } \partial \Omega\}, \quad
\hat{V} := H_0^1(\Omega) = \{ v \in H^1(\Omega): v = 0 \text{ on } \partial \Omega\}.
\end{equation*}

We suppose that $u(t) \in V$ for all $t \in I$. Then, the variational formulation \eqref{eq:fine_model} will be as follows: Find $u(t) \in V$ such that
\begin{equation}\label{eq:variational_formulation}
\begin{split}
m\left(\frac{\partial u}{\partial t}, v\right) + a(u, v) + g(u, v) = L(v; u) \quad &\text{for all } v \in \hat{V}, \quad t \in I,\\
m(u(0), v) = m(u_0, v) \quad &\text{for all } v \in \hat{V},
\end{split}
\end{equation}
where the bilinear and linear forms are following
\begin{equation}\label{eq:bilinear_forms}
\begin{gathered}
m(u, v) = \int_\Omega u v dx, \quad 
a(u, v) = \int_\Omega \epsilon \kappa(x) \nabla u \cdot \nabla v dx, \\
g(u, v) = \int_\Omega \beta \cdot \nabla u dx, \quad
L(v; u) = \int_\Omega f(u) v dx.
\end{gathered}
\end{equation}

In order to discretize \eqref{eq:variational_formulation}, we introduce the fine grid $\mathcal{T}^h$, which partitions $\Omega$ into tetrahedral cells $K^h$. We assume that $\mathcal{T}^h$ is fine enough to resolve all the heterogeneities of $\kappa(x)$. Then, we define the following finite-element spaces 
\begin{equation*}
V^h := \{ v \in V: v|_{K^h} \in P_1(K^h) \quad \forall K^h \in \mathcal{T}^h\}, \quad
\hat{V}^h := \{ v \in \hat{V}: v|_{K^h} \in P_1(K^h) \quad \forall K^h \in \mathcal{T}^h\},
\end{equation*}
where $P_1(K^h)$ is the space of all linear functions over the element $K^h$.

Therefore, we have the following Galerkin formulation of \eqref{eq:variational_formulation}: Find $u(t) \in V^h$ such that
\begin{equation}\label{eq:galerkin_formulation}
\begin{split}
m\left(\frac{\partial u}{\partial t}, v\right) + a(u, v) + g(u, v) = L(v; u) \quad &\text{for all } v \in \hat{V}^h, \quad t \in I,\\
m(u(0), v) = m(u_0, v) \quad &\text{for all } v \in \hat{V}^h,
\end{split}
\end{equation}
where $u = \sum_{j = 1}^{N_v^h} u_j \phi_j$, $u_j$ are nodal values of the solution field, $\phi_j$ are the standard piecewise linear basis functions, and $N_v^h$ is the number of the fine-grid nodes.

Moreover, one can represent \eqref{eq:galerkin_formulation} in the following matrix-vector form: 

\noindent Find $u(t) = [u_1(t), u_2(t), ..., u_{N_v^h}(t)]^T$ such that
\begin{equation}\label{eq:fine_system_continuous}
\begin{split}
M \frac{du}{dt} + (A + G) u &= F(u), \quad \text{for all } t \in I,\\
M u(0) &= \hat{u}_0, 
\end{split}
\end{equation}
where the matrices and vectors are defined as follows
\begin{equation}\label{eq:fine_matrices_and_vectors}
\begin{gathered}
M = [m_{ij}], \quad m_{ij} = m(\phi_j, \phi_i), \quad
A = [a_{ij}], \quad a_{ij} = a(\phi_j, \phi_i), \quad
G = [g_{ij}], \quad g_{ij} = g(\phi_j, \phi_i), \\
F = [f_i], \quad f_i = L(\phi_i; u), \quad
\hat{u}_0 = [\hat{u}_{0, i}], \quad
\hat{u}_{0, i} = m(u_0, \phi_i).
\end{gathered}
\end{equation}
\subsection{Backward Euler method}
Let us apply the backward Euler method for temporal approximation with $N_t$ uniform time steps of size $\tau$ such that $T = \tau N_t$. In this way, we have the following fully discrete problem: Find $u^n = [u_1^n, u_2^n, ..., u_{N_v^h}^n]^T$ such that
\begin{equation}\label{eq:fine_system_discrete}
\begin{split}
M \frac{u^{n} - u^{n - 1}}{\tau} + (A + G) u^n &= F^{n - 1}, \quad n = 1, 2, ..., N_t,\\
M u^0 &= \hat{u}_0, 
\end{split}
\end{equation}
where $u^{n} = u(t_n)$, $F^{n - 1} = F(u(t_{n - 1}))$, and $t_n = \tau n$.

\subsection{Exponential integration}\label{sec:eifinemesh}

Let us reformulate \eqref{eq:fine_system_continuous} in the following way
\begin{equation}\label{eq:cont_time_matrix_form_2}
\begin{split}
\frac{d u}{d t} + M^{-1} A u &= M^{-1} (F(u) - G u),\\
u(0) &= u_0,
\end{split}
\end{equation}
where $u_0 = M^{-1} \hat{u}_0$.

Note that this matrix ODE system has the exact solution called the variation of constants
\begin{equation}\label{eq:variation_of_constants}
u (t^n) = e^{-\tau M^{-1} A} u(t^{n - 1}) + \int_{0}^{\tau} e^{(l - \tau) M^{-1} A} M^{-1}\left[F(u(t^{n - 1} + l)) - G u(t^{n - 1} + l)\right] dl.
\end{equation}

We can approximate the integral term by the exponential quadrature rule with $s$ quadrature points $c_i \in [0, 1]$ in the following way
\begin{equation}\label{eq:exp_quadrature_rule}
u^n = e^{-\tau M^{-1} A} u^{n - 1} + \tau \sum_{i = 1}^{s} b_i (-\tau M^{-1} A) M^{-1} ( F_i - G u_i ),
\end{equation}
where $u^n \approx u(t^n)$, $F_i = F(u(t^{n - 1} + c_i \tau))$, $G u_i = G u (t^{n - 1} + c_i \tau)$, and
\begin{equation}
\nonumber
b_i (-\tau M^{-1} A) = \int_{0}^{1} e^{-\tau (1 - \theta) M^{-1} A} l_i(\theta) d\theta.
\end{equation}
Note that $l_i(\theta)$ are the Lagrange interpolation polynomials.

The coefficients $b_i (z)$ can also be considered as linear combinations of $\varphi$-functions defined as
\begin{equation}\label{eq:varphi_definition}
\begin{split}
\varphi_0 (z) = e^z, \quad
\varphi_{p + 1} (z) = \frac{1}{z} \left( \varphi_p(z) - \frac{1}{p!} \right).
\end{split}
\end{equation}
with the following recurrence relation
\begin{equation}\label{eq:varphi_definition_2}
\varphi_p (z) = \int_0^1 e^{(1 - \theta) z} \frac{\theta^{p - 1}}{(p - 1)!} d\theta, \quad p \geq 1.
\end{equation}

Let us take $s = 1$ and $c_1 = 0$. Then, we have $b_1(-\tau M^{-1} A) = \int_0^1 e^{-\tau (1 -\theta) M^{-1} A} d \theta = \varphi_1(-\tau M^{-1} A)$, according to \eqref{eq:varphi_definition_2}. By taking $e^z = z \varphi_1(z) + 1$, $F_1 \approx F^{n - 1} = F(u^{n - 1})$, and $G u_1 \approx G u^{n - 1}$ in \eqref{eq:exp_quadrature_rule}, we get
\begin{equation}\label{eq:exp_euler_scheme}
u^n = u^{n - 1} + \tau \varphi_1 (-\tau M^{-1} A) M^{-1}([F^{n - 1} - G u^{n - 1}] - A u^{n - 1}),
\end{equation}

Next, by the symmetry and positive definiteness of $A$ and $M$, we have \cite{contreras2023exponential}
\begin{equation}
\nonumber
-\tau M^{-1} A = Q D Q^T M
\end{equation}
and
\begin{equation}
\nonumber
\varphi_1(-\tau M^{-1} A) = Q \varphi_1 (D) Q^T M,
\end{equation}
where $Q$ is formed with eigenvectors of $-\tau A q = \lambda M q$ arranged as columns.

Therefore, we obtain
\begin{equation}\label{eq:exp_euler_scheme_2}
u^n = u^{n - 1} + \tau Q \varphi_1(D) Q^T (F^{n - 1} - (A + G) u^{n - 1}).
\end{equation}

In the next section, we present a meshfree generalized multiscale finite element method combined with the exponential integration.

\section{Meshfree generalized multiscale finite element method}\label{sec:meshfreegmsfem}
In this section, we present the construction of a reduced-order approximation based on a meshfree coarse representation. 
We briefly describe the meshfree generalized multiscale finite element method (MFGMsFEM), where coarse basis functions are constructed from local spectral problems using only a point cloud of coarse nodes. This approach avoids the need for explicit coarse meshes while retaining the key features of classical GMsFEM. For more details we refer to \cite{nikiforov2025meshfree} and references therein. 

\subsection{Coarse scale}

First, we introduce the coarse-scale points $\{x_i\}_{i=1}^{N_v^H}$, where $N_v^H$ is the number of coarse points. We proceed as in 
\cite{djulustan2024meshfree, nikiforov2023meshfree,nikiforov2023modeling, nikiforov2023meshfreeRichards} and references therein. We partition the computational domain $\Omega$ into overlapping coarse elements $S_i$ constructed around each $x_i^H$ so that $\Omega \subset \cup_{i = 1}^{N_v^H} S_i$. In addition, we assume that each $S_i$ can be represented as a union of fine-grid elements. More precisely, we have
\begin{equation}
S_i = \cup \left\{ K^h \in \mathcal{T}^h: \text{max}_{y \in K^h} \| y - x_i \| \leq r_i \right\},
\end{equation}
where $r_i$ denotes the radius of the coarse element $S_i$.

In the meshfree generalized multiscale finite element method, we construct multiscale basis functions $\phi_{ik}^{\text{ms}}$, supported in $S_i$, which capture the fine-scale heterogeneities. Using these basis functions, the solution can be approximated as follows
\begin{equation}
u_{\text{ms}}(t, x) = \sum_{i,k} u_{\text{ms}, ik} (t) \phi_{ik}^{\text{ms}} (x),
\end{equation}
where the subscripts $i$ and $k$ correspond to the coarse point and basis function indices, respectively.

\begin{figure}[!htb]
\begin{center}
\begin{minipage}[h]{0.49\linewidth}
\center{\includegraphics[width=\linewidth]{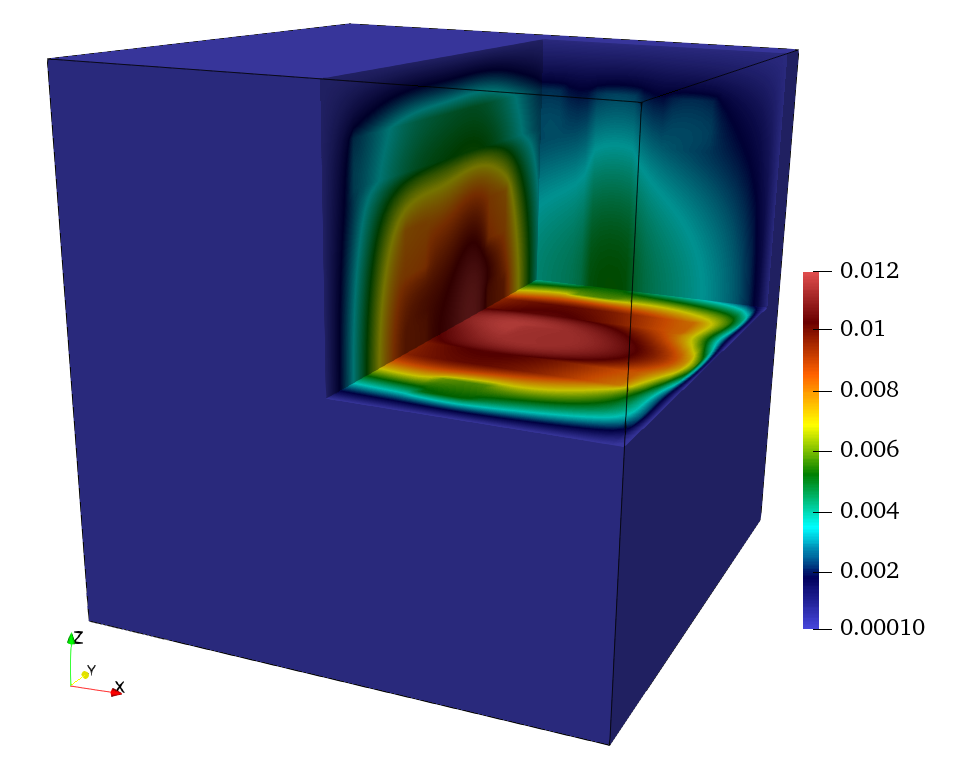}}
\end{minipage}
\begin{minipage}[h]{0.49\linewidth}
\center{\includegraphics[width=\linewidth]{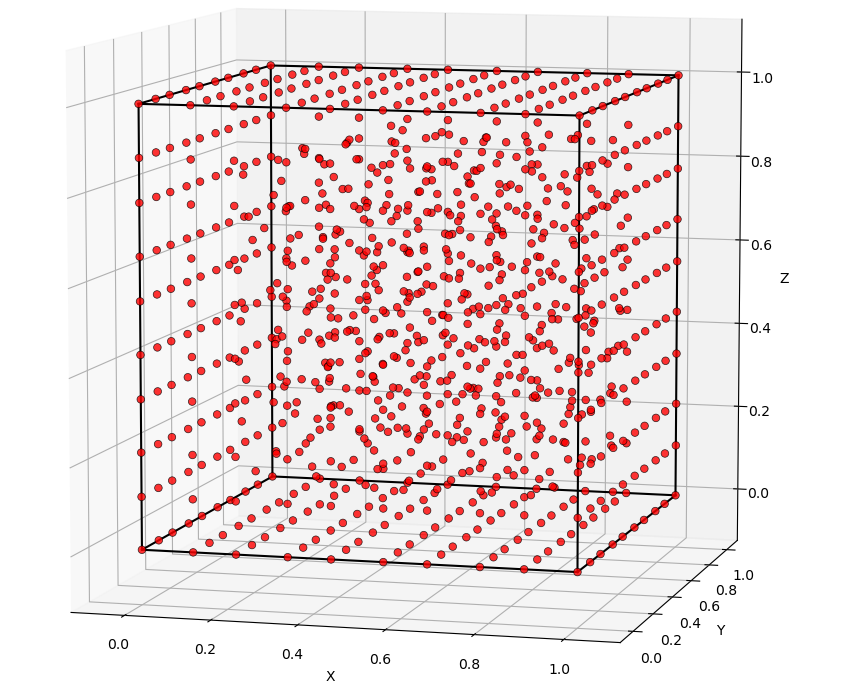}}
\end{minipage}
\end{center}
\caption{Illustration of the density distribution function $\rho(x)$ (left) and the coarse-scale points $\{x_i\}_{i=1}^{N_v^H}$ (right).}
\label{fig:point_cloud}
\end{figure}

To construct the distribution of coarse-scale points $\{x_i\}_{i=1}^{N_v^H}$ (see Fig. \ref{fig:point_cloud}) and the associated radii $r_i$, we employ the probabilistic algorithms proposed in \cite{du2002meshfree}. The parameters used in our computations are: number of iterations $N{\text{iter}} = 200$, number of random samples per iteration $q = 50000$, update coefficients $\alpha_1 = \alpha_2 = \beta_1 = \beta_2 = 0.5$, uniform pseudo‑point spacing $h_i = 0.005$, and overlap factor $\gamma = 2$. The algorithms are summarized in Algorithms~\ref{alg:cvt} and~\ref{alg:radii}. We compute the density distribution function $\rho(x)$ used in the CVT construction by solving the following problem
\begin{equation}
\begin{split}
\alpha (-\epsilon \div (\kappa(x) \nabla \rho) + \beta \cdot \nabla \rho) + \rho = u_0(x) \quad &\text{in } \Omega,\\
\rho = \varepsilon \quad &\text{on } \partial \Omega,
\end{split}
\end{equation}
where $\alpha$ is the smoothing parameter and $\varepsilon$ is a small positive number. Here, we set the smoothing parameter $\alpha = 10^{-1}$ and the boundary value $\varepsilon = 10^{-4}$; these values are chosen based on numerical experiments to achieve a balance between regularity and accuracy. The above problem is used to smooth a initial distribution while accounting for permeability $\kappa(x)$ and convective effects $\beta$, thereby obtaining a density function $\rho(x)$ that better reflects the underlying physical features for the subsequent centroidal Voronoi tessellation construction.
See \cite{nikiforov2025meshfree} and references therein for more details on the constructions. 

\subsection{Local multiscale basis functions}

In our previous work \cite{nikiforov2025meshfree}, multiscale basis functions were constructed within the meshfree GMsFEM framework using local spectral problems tailored to diffusion-dominated settings. In the present work, we extend this construction to advection--diffusion problems, where the presence of transport effects requires additional considerations in the design of the coarse space. In particular, we introduce two types of multiscale basis functions that incorporate advection either at the snapshot level or directly within the spectral problem, allowing for improved approximation properties in advection-influenced regimes. 

Here we use ideas proposed in \cite{xie2025time}. In particular,  we consider two types of multiscale basis functions. In the first type, we construct the basis functions using the spectral problem that does not consider advection effects. The advection is considered only in the construction of snapshot functions. In contrast, the second type of basis functions employs an alternative spectral problem to incorporate the advection effects.

\paragraph{Type 1 multiscale basis functions.}

First, let us consider the construction of snapshot functions designed to account for advection. For this purpose, we solve the following local problems. For $i = 1, 2, ..., N_v^H$ and all $j \in J(S_i)$, find $\phi_j^{\text{snap}, S_i}$ such that
\begin{equation}\label{eq:snap1}
\begin{split}
- \epsilon \div \left( \kappa(x) \nabla \phi_j^{\text{snap}, S_i} \right) + \beta \cdot \nabla \phi_j^{\text{snap}, S_i} &= 0 \quad \text{in } S_i,\\
\phi_j^{\text{snap}, S_i} &= \delta_j(x) \quad \text{on } \partial S_i,
\end{split}
\end{equation}
where $\delta_j(x) = \delta_{jk}$ for all $k \in J(S_i)$, and $J(S_i)$ is a set of all the fine-grid nodes on $\partial S_i$.

After computing the snapshot functions, we construct the following snapshot spaces and projection matrices
\begin{equation}\label{eq:snapshot_space}
\begin{split}
V_{\text{snap}}(S_i) = \text{span} \{ \phi_j^{\text{snap}, S_i}: ~1 \leq j \leq J_i\},\\
R_i^{\text{snap}} = [\phi_1^{\text{snap}, S_i}, \phi_2^{\text{snap}, S_i}, ..., \phi_{J_i}^{\text{snap}, S_i}]^T,
\end{split}
\end{equation}
where $J_i$ is the number of the fine-grid nodes on $\partial S_i$.

Then, we solve the following local spectral problems in the snapshot spaces. For each $i = 1,2, ..., N_v^H$, find $\lambda_k^{S_i}$ and $\tilde{\phi}_k^{S_i}$ such that
\begin{equation}\label{eq:spectral_problem_1}
\tilde{A}_i \tilde{\phi}_k^{S_i} = \lambda_k^{S_i} \tilde{B}_i \tilde{\phi}_k^{S_i},
\end{equation}
where the matrices are given as follows
\begin{equation*}
\begin{gathered}
\tilde{A}_i = R_i^{\text{snap}} A_i (R_i^{\text{snap}})^T, \quad 
\tilde{B}_i = R_i^{\text{snap}} B_i (R_i^{\text{snap}})^T, \\
A_i = [a_{i, kl}], \quad a_{i, kl} = \int_{S_i} \epsilon \kappa(x) \nabla \phi_l \cdot \nabla \phi_k dx, \quad
B_i = [b_{i, kl}], \quad b_{i, kl} = \int_{S_i} \epsilon \kappa(x) \phi_l \phi_k dx,
\end{gathered}
\end{equation*}

We sort the eigenvalues in ascending order and take the first $N_{\text{b}}$ eigenvectors as basis functions. Using the snapshot projection matrices, we project them from the snapshot spaces $\phi_k^{S_i} = (R_i^{\text{snap}})^T \tilde{\phi}_k^{S_i}$.

\paragraph{Type 2 multiscale basis functions.}

In the second type, we construct the snapshot spaces and corresponding projection matrices in the same way. However, we consider the following alternative spectral problems accounting for advection effects. For $i = 1,2,...,N_v^H$, find $\lambda_k^{S_i}$ and $\phi_k^{S_i}$ such that
\begin{equation}\label{eq:spectral_problem_2}
\tilde{D}_i \phi_k^{S_i} = \lambda_k^{S_i} \tilde{N}_i \phi_k^{S_i}, \\
\end{equation}
where the matrices defined as follows
\begin{equation*}
\begin{gathered}
\tilde{D}_i = R_i^{\text{snap}} D_i (R_i^{\text{snap}})^T, \quad
\tilde{N}_i = R_i^{\text{snap}} N_i (R_i^{\text{snap}})^T, \\
D_i = (A_i + G_i)^T M_i (A_i + G_i), \quad
N_i = M_i^2,\\
A_i = [a_{i, kl}], \quad a_{i, kl} = \int_{S_i} \epsilon \kappa(x) \nabla \phi_l \cdot \nabla \phi_k dx, \quad
G_i = [g_{i, kl}], \quad g_{i, kl} = \int_{S_i} \beta \cdot \nabla \phi_l \phi_k dx, \\
M_i = [m_{i, kl}], \quad m_{i, kl} = \int_{S_i} \phi_l \phi_k  dx ,
\end{gathered}
\end{equation*}

Again, we sort the eigenvalues in ascending order and choose the first $N_{\text{b}}$ eigenvectors as basis functions.

\subsection{Global formulation}

We use the following shape functions to ensure conformality of the multiscale basis functions
\begin{equation}
W_i (x) = \frac{\eta_i(x)}{\sum_{j = 1}^N \eta_j(x)},
\end{equation}
where $\eta_j(x)$ are the kernel functions
\begin{equation}
\eta_i(x) = 2 \begin{cases}
2/3 + 4 (r - 1) r^2, & r \leq 0.5,\\
4/3 (1 - r)^3, & 0.5 \leq r \leq 1,\\
0, & 1 \leq r.
\end{cases}
\end{equation}
Here, $r = \frac{\| x - x_k \|}{r_k}$ is the normalized distance.

\begin{figure}[!htb]
\begin{center}
\begin{minipage}[h]{0.49\linewidth}
\center{\includegraphics[width=\linewidth]{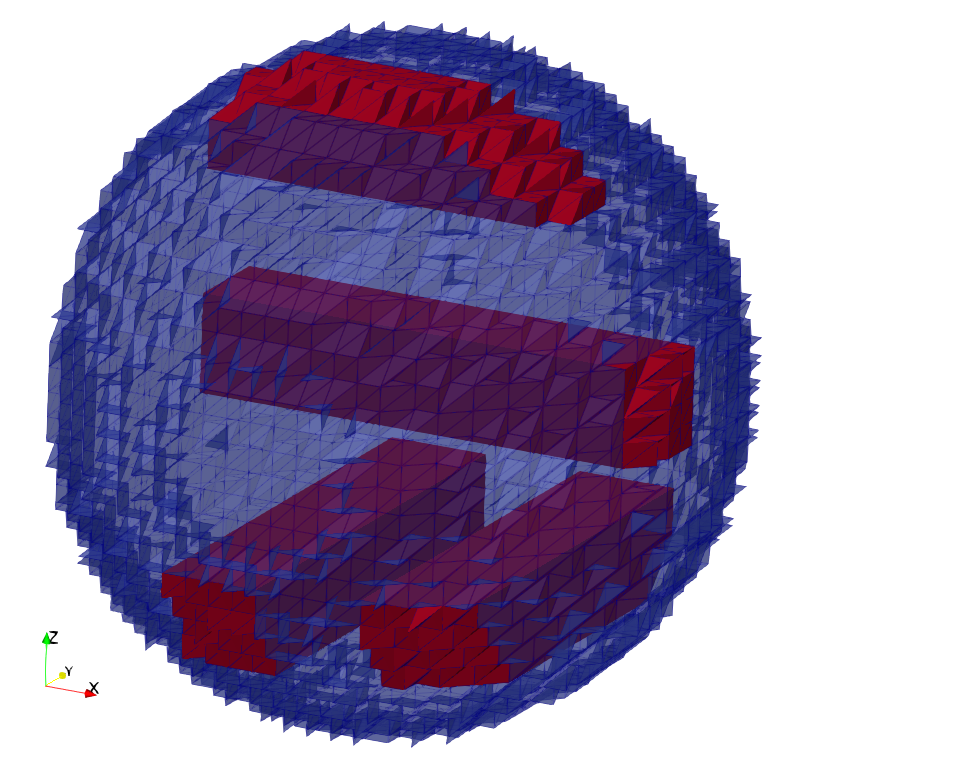}}
\end{minipage}
\begin{minipage}[h]{0.49\linewidth}
\center{\includegraphics[width=\linewidth]{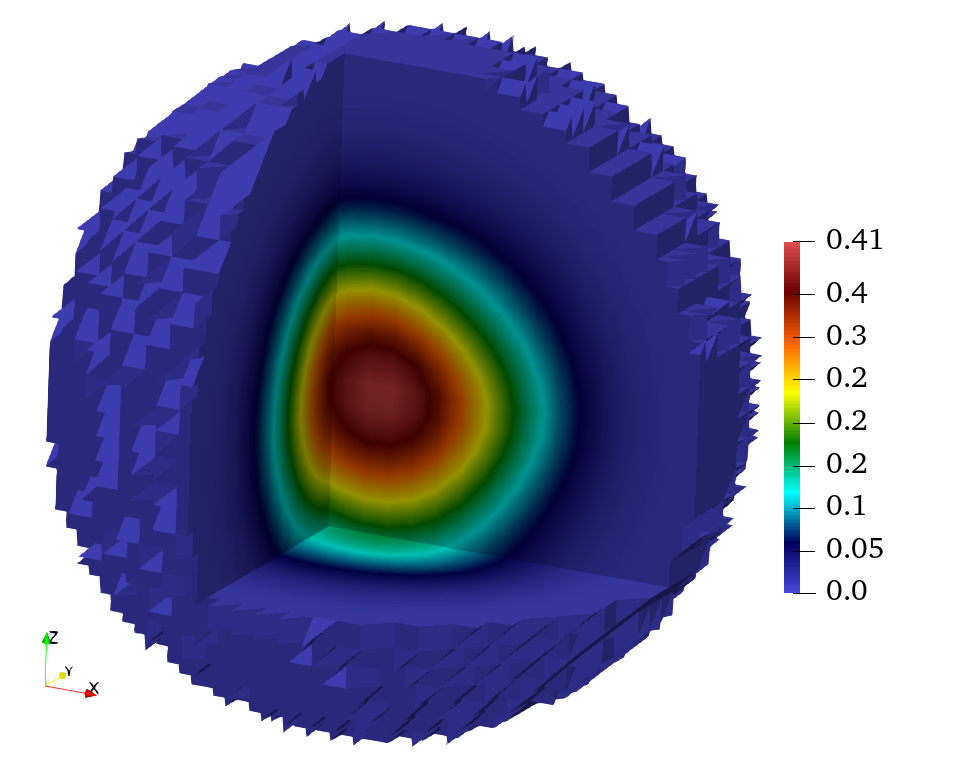}}
\end{minipage}\\
\begin{minipage}[h]{0.49\linewidth}
\center{\includegraphics[width=\linewidth]{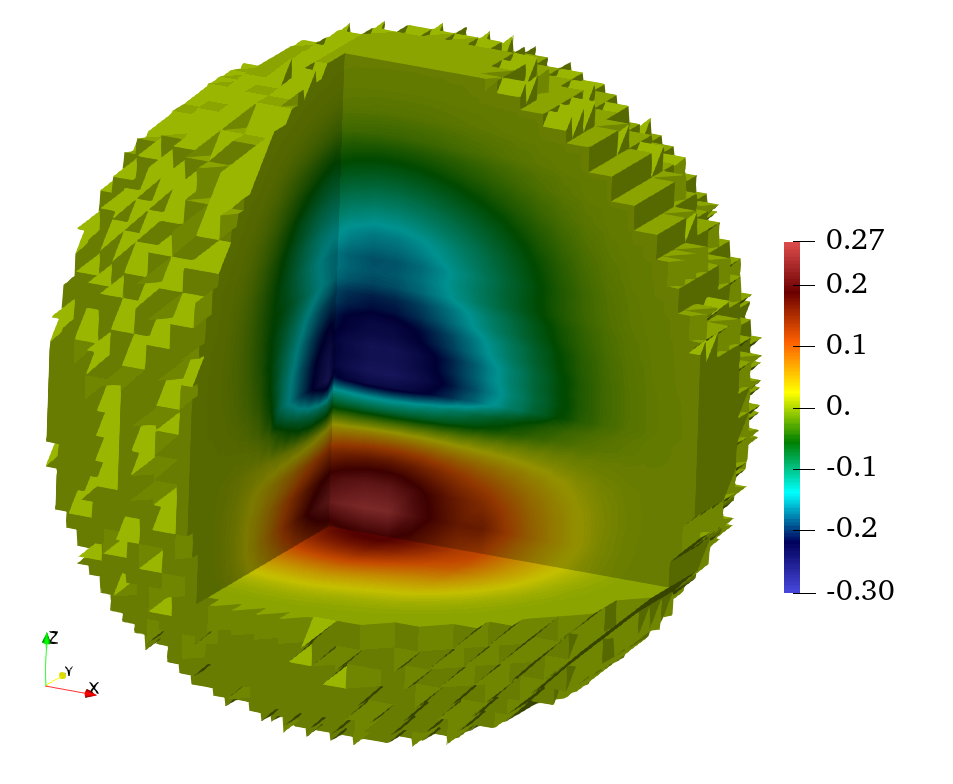}}
\end{minipage}
\begin{minipage}[h]{0.49\linewidth}
\center{\includegraphics[width=\linewidth]{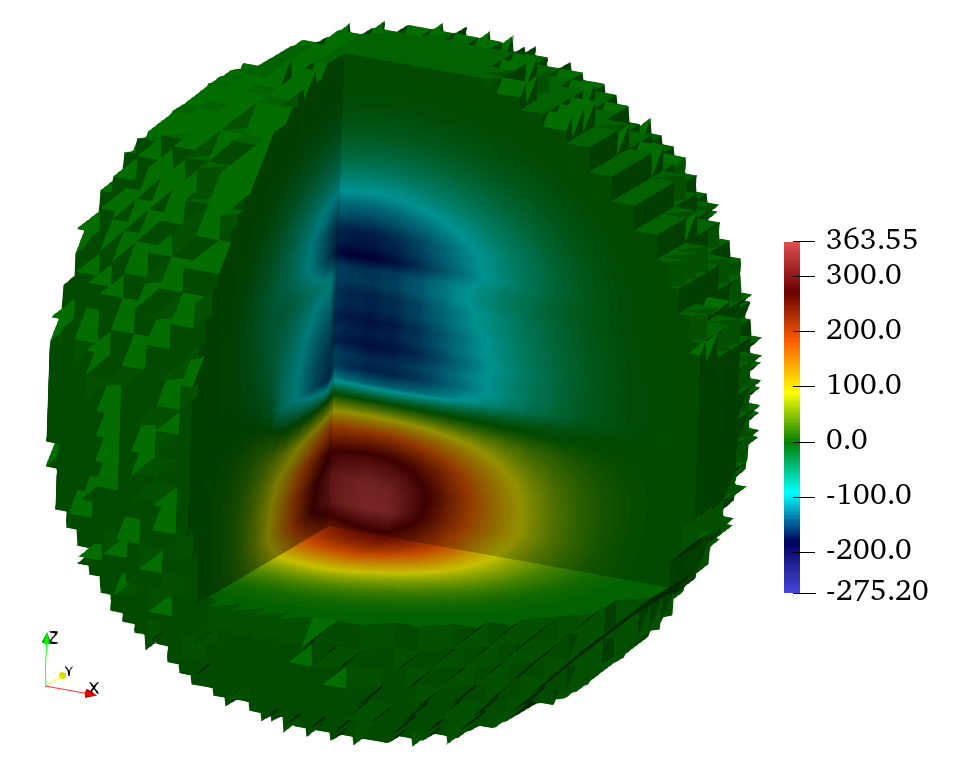}}
\end{minipage}
\end{center}
\caption{Illustration of the high-contrast coefficient $\kappa(x)$ (top left), the shape function $W_i (x)$ (top right), the type 1 basis function (bottom left), and the type 2 basis function (bottom right).}
\label{fig:basis_illustration}
\end{figure}

We multiply the eigenfunctions by the shape functions, and obtain the multiscale basis functions
\begin{equation}
\phi_{ik}^{\text{ms}} = W_i \phi_{k}^{S_i}.
\end{equation}
Examples of the resulting multiscale basis functions together with the high‑contrast coefficient and a representative shape function are shown in Fig.~\ref{fig:basis_illustration}.

Finally, we can construct the multiscale space and the projection matrix as follows
\begin{equation}
\begin{gathered}
V_{\text{ms}} = \text{span} \{ \phi_{ik}^{\text{ms}}: 1 \leq i \leq N_v^H, \, 1 \leq k \leq N_{\text{b}} \},\\
R_{\text{ms}} = [\phi_{11}^{\text{ms}}, \phi_{12}^{\text{ms}}, ..., \phi_{1N_{\text{b}}}^{\text{ms}}, \phi_{21}^{\text{ms}}, \phi_{22}^{\text{ms}}, ..., \phi_{2N_{\text{b}}}^{\text{ms}}, ...,  \phi_{N_v^H 1}^{\text{ms}}, \phi_{N_v^H 2}^{\text{ms}}, ..., \phi_{N_v^H N_{\text{b}}}^{\text{ms}}]^T.
\end{gathered}
\end{equation}

Note that the Galerkin formulation of our problem in the multiscale space can be written as
\begin{equation}\label{eq:ms_galerkin_formulation}
m\left(\frac{\partial u_{\text{ms}}}{\partial t}, v_{\text{ms}}\right) + a(u_{\text{ms}}, v_{\text{ms}}) + g(u_{\text{ms}}, v_{\text{ms}}) = l(v_{\text{ms}}; u_{\text{ms}}) \quad \text{for all } v_{\text{ms}} \in \hat{V}_{\text{ms}},\\
\end{equation}
which also has the following matrix-vector form
\begin{equation}
M_H \frac{du_{\text{ms}}}{dt} + (A_H + G_H) u_{\text{ms}} = F_H(u_{\text{ms}}),
\end{equation}
where
\begin{equation*}
\begin{gathered}
M_H = R_{\text{ms}} M (R_{\text{ms}})^T, \quad
A_H = R_{\text{ms}} A (R_{\text{ms}})^T, \quad
G_H = R_{\text{ms}} G (R_{\text{ms}})^T, \quad
F_H = R_{\text{ms}} F(u_{\text{ms}}).
\end{gathered}
\end{equation*}

A standard approach for temporal discretization is to use the backward Euler method.
\begin{equation}
M_H \frac{u_{\text{ms}}^n - u_{\text{ms}}^{n - 1}}{\tau} + (A_H + G_H) u_{\text{ms}}^n = F_H(u_{\text{ms}}^{n - 1}).
\end{equation}

However, as discussed above and in 
many previous works, e.g., \cite{nikiforov2025meshfree, contreras2023exponential, xie2025time, xie2025multiscale}, 
such a time discretization scheme may have stability issues for problems in multiscale media. This is the main reason to introduct the exponential integration approach but a explained in before and in previous works
(\cite{nikiforov2025meshfree, contreras2023exponential, xie2025time, xie2025multiscale} this is impractical for high-contras problems applications, even more for 3D problems.  See Section \ref{sec:eifinemesh} and problem \ref{eq:exp_euler_scheme_2}. 
To reduce the computational cost, we approximate the matrix function by upscaling it with multiscale projection matrix
\begin{equation}\label{eq:mgmsfem_ei_2}
u^n = u^{n - 1} + \tau R_{\text{ms}}^T Q_H \varphi_1 (D_H) Q_H^T R_{\text{ms}} (F^{n - 1} - (A + G) u^{n - 1}),
\end{equation}
where $Q_H$ is formed with the coarse-scale eigenvectors (as its columns) of the following coarse-scale spectral problem
\begin{equation}
-\tau A_H q_H = \lambda M_H q_H.
\end{equation}

We use the orthogonal projection for the initial condition
\begin{equation}
\hat{u}_H = R_{\text{ms}}^T M_H^{-1} R_{\text{ms}} M u_H.
\end{equation}

For more details, deductions and motivation we refer to \cite{nikiforov2025meshfree, xie2025time, contreras2023exponential}.

\begin{remark}
The theoretical analysis of the proposed methodology is not revisited in this work, as it is already covered in our previous paper \cite{nikiforov2025meshfree}, where the stability and convergence properties of the meshfree GMsFEM combined with exponential integration were established. Additional theoretical developments for exponential integration in multiscale settings can be found in \cite{contreras2023exponential, xie2025time, xie2025multiscale}, which build upon the framework introduced in \cite{poveda2024second,poveda2024edge}. In the present work, we focus on the computational extension of the method to three-dimensional problems and on the design of multiscale basis functions tailored to advection-dominated regimes.
\end{remark}

We now proceed to assess the performance of the proposed meshfree multiscale exponential integration approach through a set of three-dimensional numerical experiments, focusing on accuracy, stability, and computational efficiency.

\section{Numerical results}\label{sec:numerical_results}

In this section, we present a set of three-dimensional numerical experiments to evaluate the performance of the proposed meshfree generalized multiscale exponential integration (MFGMsFEM-EI) method. The experiments are designed to assess the accuracy, stability, and computational efficiency of the method in the presence of heterogeneous and high-contrast diffusion coefficients, as well as nonlinear and time-dependent source terms.

We consider the following advection--diffusion problem:
\begin{equation}
\frac{\partial u}{\partial t} - \epsilon \, \text{div}(\kappa(x)\nabla u) + \beta \cdot \nabla u = f \quad \text{in } \Omega \times I,
\end{equation}
with homogeneous Dirichlet boundary conditions and initial condition
\[
u(0,x) = x_1(1-x_1)x_2(1-x_2)x_3(1-x_3),
\]
where $\Omega = (0,1)^3$ and the velocity field is given by
\[
\beta = [2+\sin(4\sqrt{2}\pi x_2), 0, 0].
\]

The heterogeneous diffusion coefficient $\kappa(x)$ exhibits a channelized structure, as shown in Figure~1, and is characterized by different contrast levels. The fine-grid discretization uses a structured mesh of size $50 \times 50 \times 50$, which is sufficiently fine to resolve all heterogeneities. The coarse-scale discretization for MFGMsFEM employs $10 \times 10 \times 10$ coarse regions. The final simulation time is $T = 0.2$, and we use 50 time steps for the multiscale method and 20000 time steps for the reference fine-grid solution.

We consider two types of source terms:
\[
f_1 = 100x_1(0.1 - x_1)x_2(0.1 - x_2)x_3(0.1 - x_3)\sin(t), \quad
f_2 = u(u - 1)(u + 1),
\]
leading to both linear and nonlinear problem settings.

To thoroughly assess the method, we define eight test cases combining different diffusion contrasts, advection-diffusion regimes, and source terms:
\begin{itemize}
    \item Examples 1--4 correspond to $\epsilon = 1$ with contrasts 1000 and 10.
    \item Examples 5--8 correspond to $\epsilon = 1/20$, representing more advection-dominated regimes.
    \item In each case, we consider both $f_1$ and $f_2$.
\end{itemize}

These examples are designed to explore a broad range of physical regimes, including diffusion-dominated, advection-dominated, linear, and nonlinear behaviors. While we report results for all cases in terms of error metrics, we focus the qualitative analysis on representative examples to illustrate the main features of the proposed method.

\begin{figure}[!htb]
\begin{center}
\begin{minipage}[h]{0.49\linewidth}
\center{\includegraphics[width=\linewidth]{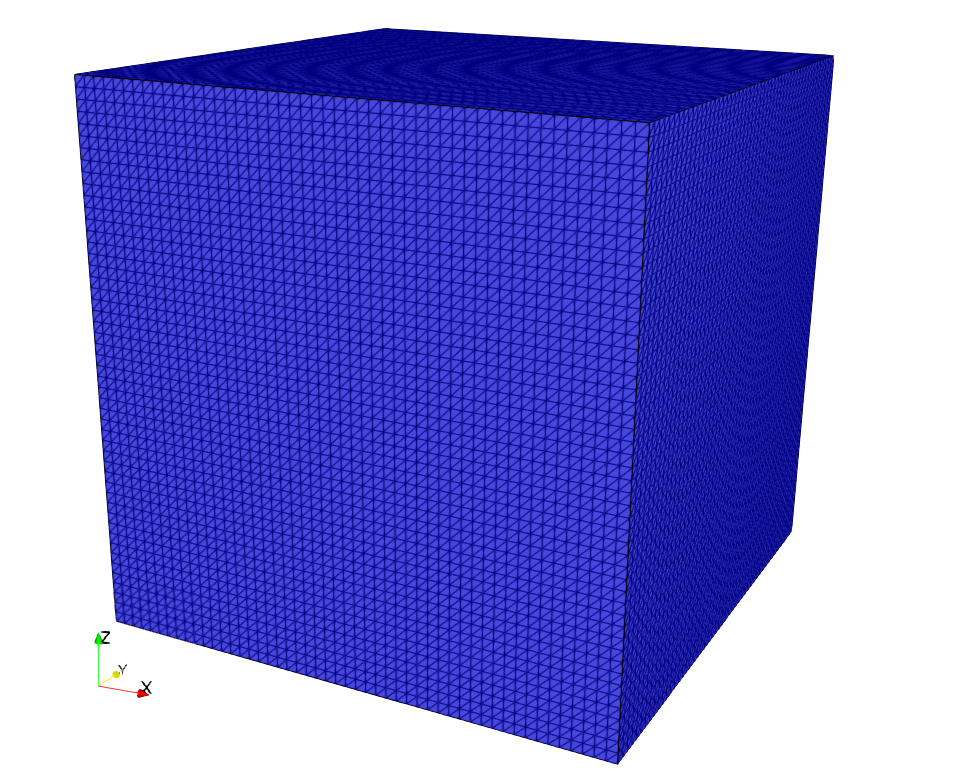}}
\end{minipage}
\begin{minipage}[h]{0.49\linewidth}
\center{\includegraphics[width=\linewidth]{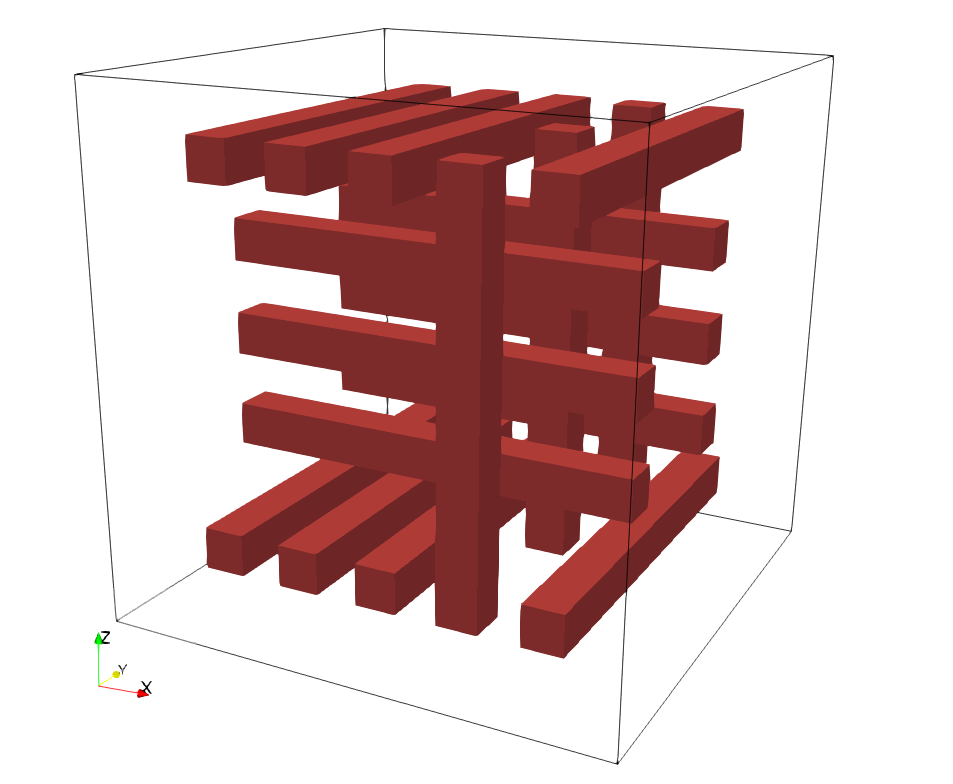}}
\end{minipage}
\end{center}
\caption{Fine grid (left) and channels (right).}
\label{fig:mesh}
\end{figure}

Figure~1 shows the fine-grid discretization (left) and the channelized permeability field (right). The complex geometry and high contrast of $\kappa(x)$ highlight the multiscale nature of the problem and justify the use of multiscale model reduction techniques.

\medskip

Figures~2 and~3 present representative solution profiles for Examples 1 and 7, respectively. In each figure, we compare the reference fine-grid solution with the multiscale approximations obtained using MFGMsFEM with standard time discretization and with exponential integration. 
\begin{figure}[!htb]
\begin{center}
\begin{minipage}[h]{0.49\linewidth}
\center{\includegraphics[width=\linewidth]{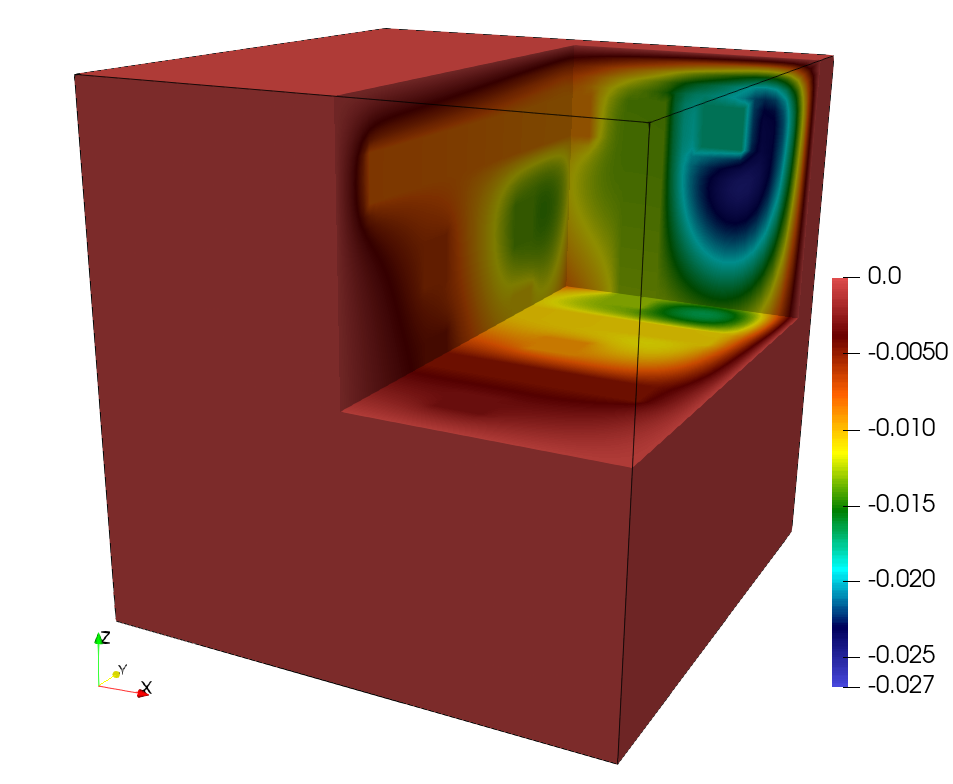}}
\end{minipage}
\begin{minipage}[h]{0.49\linewidth}
\center{\includegraphics[width=\linewidth]{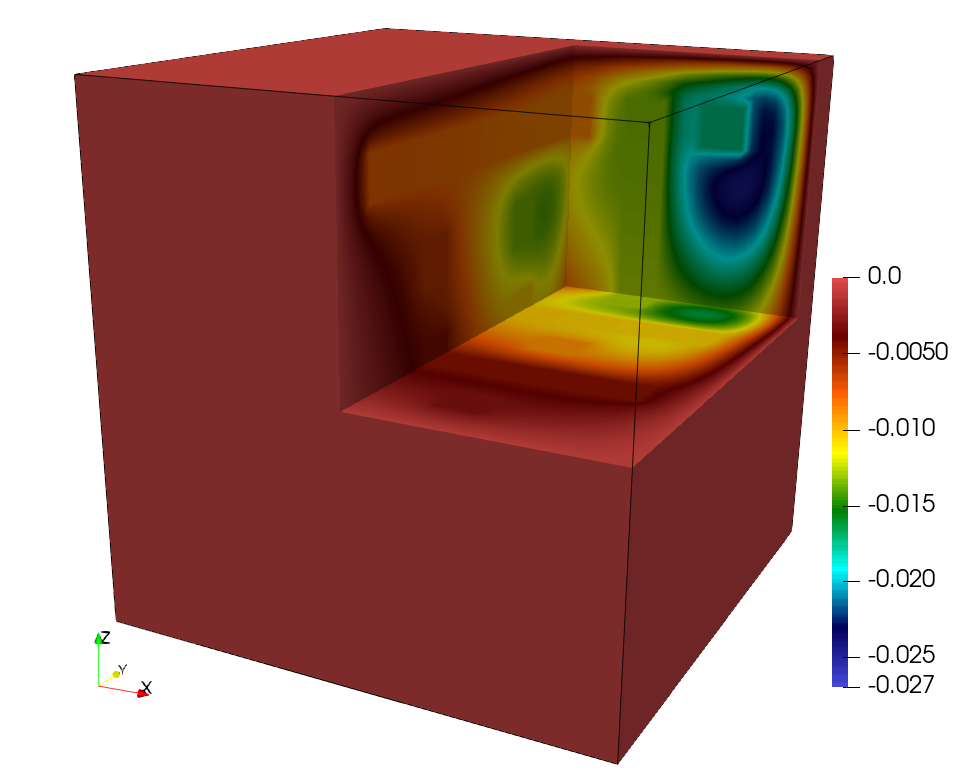}}
\end{minipage}\\
\begin{minipage}[h]{0.49\linewidth}
\center{\includegraphics[width=\linewidth]{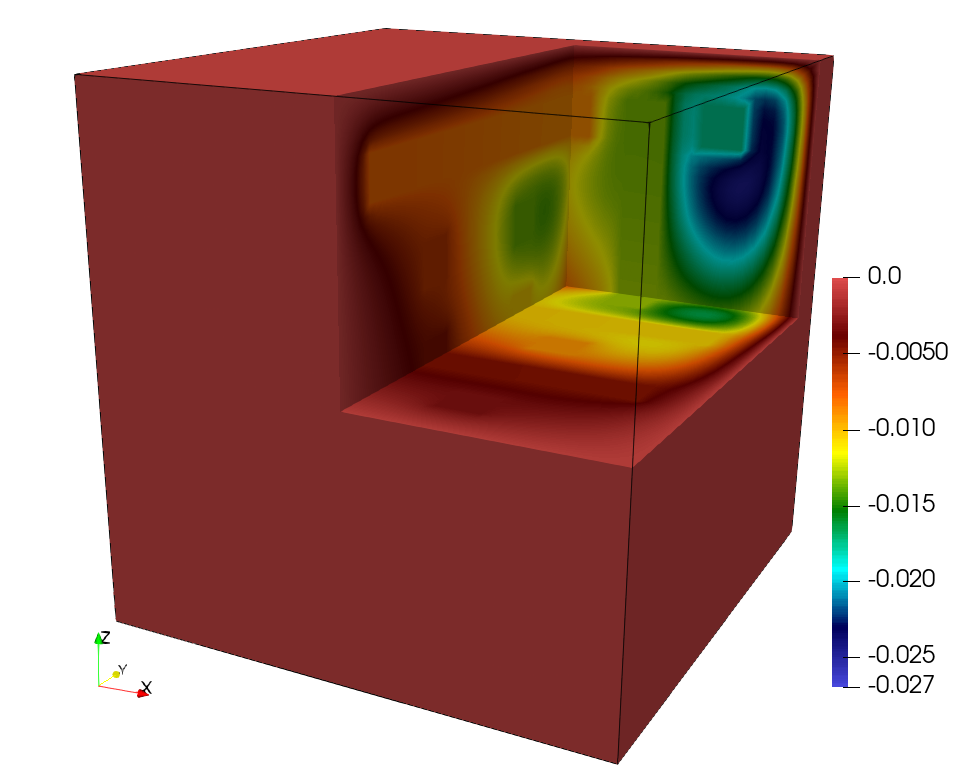}}
\end{minipage}
\end{center}
|\caption{Solution of Problem 1: fine-grid solution (Left Top), MFGMsFEM solution-FD $basis_1$ solution with 20 basis functions (Right Top) and MFGMsFEM solution-EI $basis_1$ solution with 20 basis functions (Bottom).}
\label{fig:solution1}
\end{figure}

\begin{figure}[!htb]
\begin{center}
\begin{minipage}[h]{0.49\linewidth}
\center{\includegraphics[width=\linewidth]{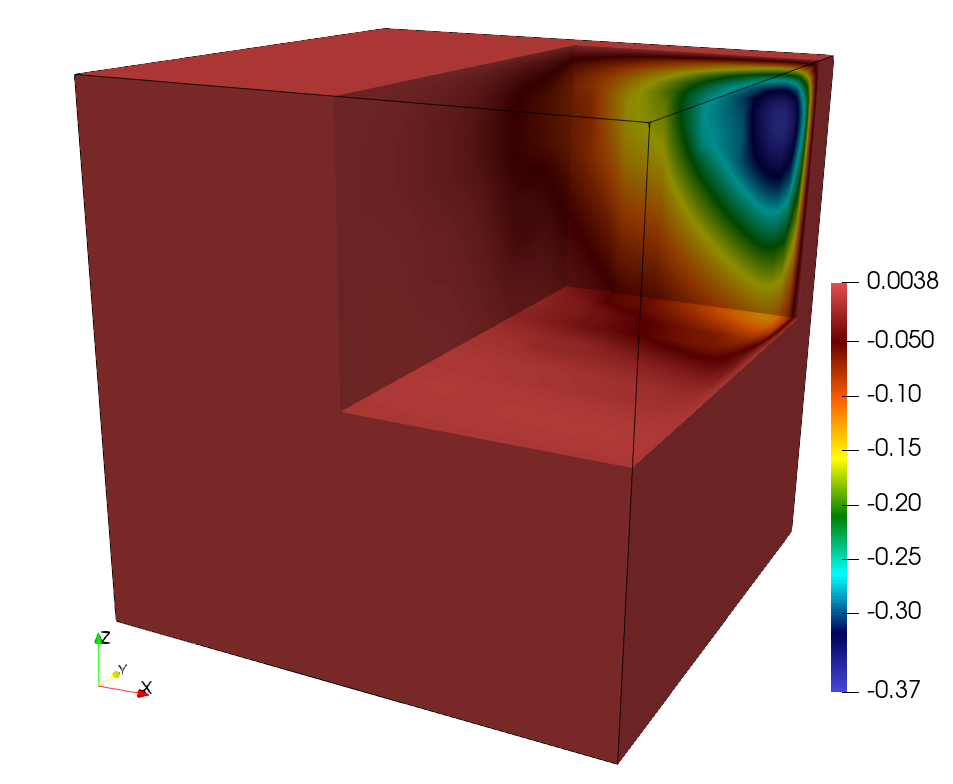}}
\end{minipage}
\begin{minipage}[h]{0.49\linewidth}
\center{\includegraphics[width=\linewidth]{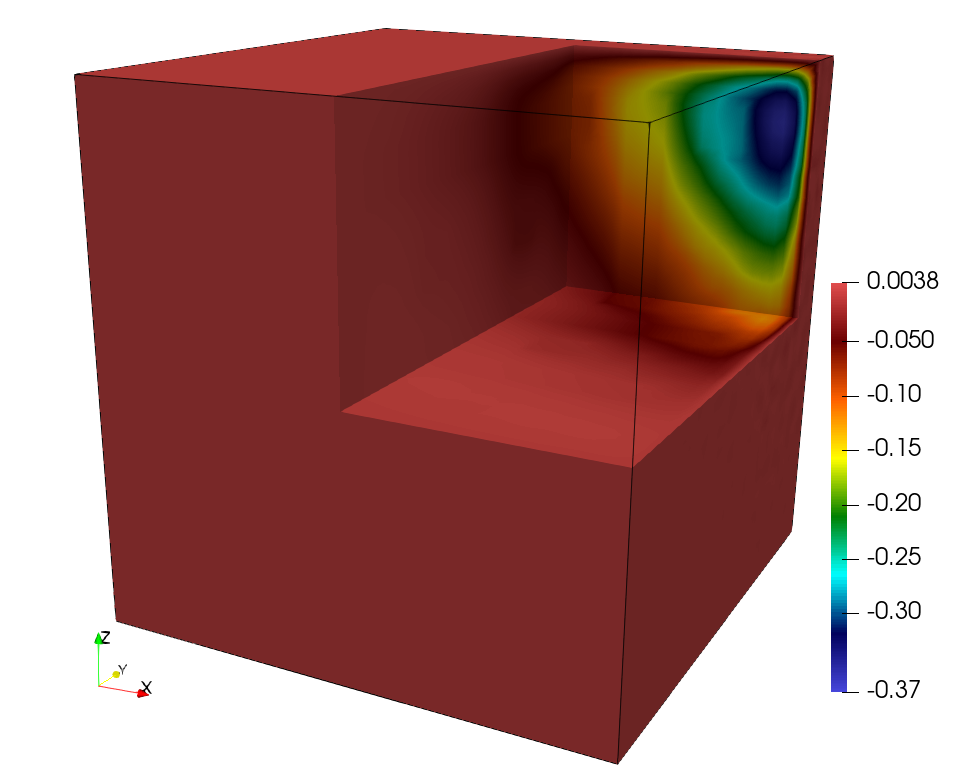}}
\end{minipage}\\
\begin{minipage}[h]{0.49\linewidth}
\center{\includegraphics[width=\linewidth]{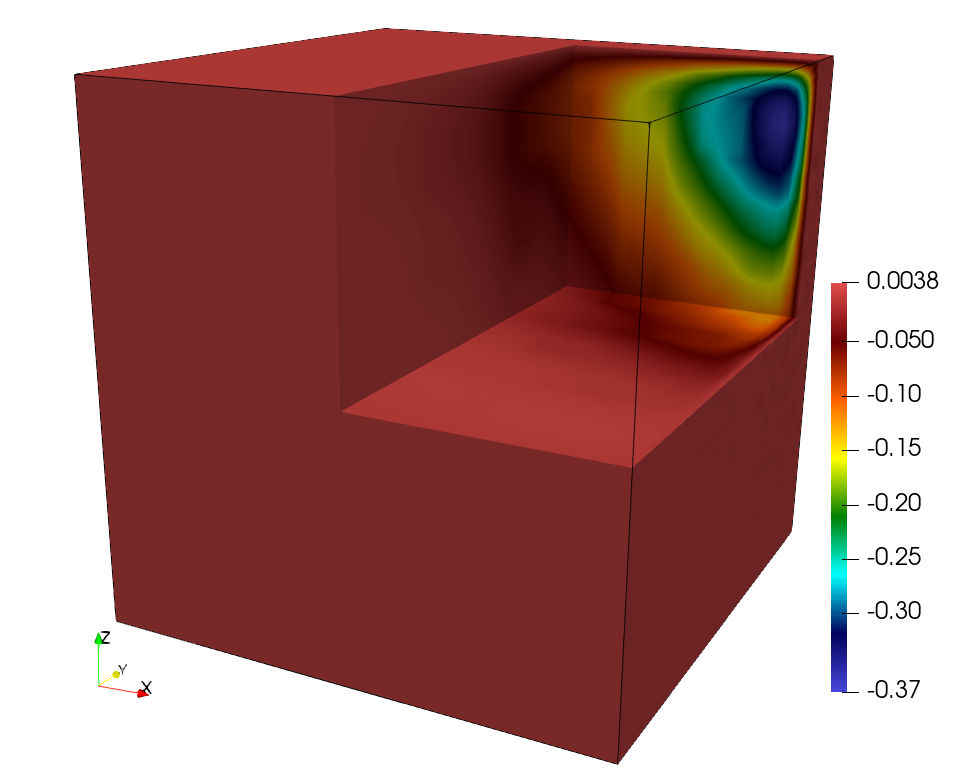}}
\end{minipage}
\end{center}
\caption{Solution of Problem 7: fine-grid solution (left top), MFGMsFEM solution-FD $basis_1$ solution with 20 basis functions (right top) and MFGMsFEM solution-EI $basis_1$ solution with 20 basis functions (bottom).}
\label{fig:solution7}
\end{figure}

In Figure~2 (Example 1), corresponding to a diffusion-dominated regime with high contrast, both multiscale approaches capture the main features of the solution. However, the exponential integration method provides a more accurate approximation, particularly in regions with sharp gradients, demonstrating its improved stability properties.

In Figure~3 (Example 7), corresponding to a more advection-dominated regime, the advantages of the exponential integration approach become more pronounced. The standard discretization shows visible discrepancies compared to the reference solution, while the MFGMsFEM-EI solution maintains good agreement, indicating better handling of transport effects and stiffness.

\medskip

In the following, we use the notation GMsFEM-FD to denote the meshfree generalized multiscale finite element method combined with a standard time discretization scheme (forward/backward Euler-type), and GMsFEM-EI to denote the proposed method where exponential integration is employed for time discretization. We compare the impact of the temporal integration strategy on the accuracy and stability of the multiscale problems, particularly in stiff and advection-dominated regimes caused by high-contrast problem.

To estimate the accuracy of the multiscale approaches, we use the following relative $L^2$ and $H^1$ errors
\[ 
{||e||}_{L^2} = \frac{{||p_1-p_2||}_{L^2}^w}{{||p_1||}_{L^2}^w}\cdot100\%, \quad
{||e||}_{H^1} = \frac{{||p_1-p_2||}_{H^1}^w}{{||p_1||}_{H^1}^w}\cdot100\%.
\]
Here, ${||p||}_{L^2}^w = \sqrt{\int_{\Omega} \kappa(x) p^2 \,dx}$, and ${||p||}_{H^1}^w = \sqrt{\int_{\Omega} \kappa(x) \nabla p \cdot \nabla p \,dx}$, where $p_1$ denotes the reference finite-element solution, and $p_2$ is a multiscale solution (MFGMsFEM-FD or MFGMsFEM-EI).

Figures~4 to 7 report the relative $L^2$ and $H^1$ errors at the final time $T = 0.2$ for all test cases. Each figure corresponds to a pair of examples with the same physical regime but different source terms.

Overall, the results show that the proposed method achieves low relative errors across all test cases, even with a significantly reduced number of degrees of freedom and larger time steps compared to the fine-grid reference solution. The exponential integration approach consistently outperforms the standard time discretization, particularly in high-contrast and advection-dominated regimes.

\begin{figure}[!htb]
\begin{center}
\begin{minipage}[h]{0.49\linewidth}
\center{\includegraphics[width=\linewidth]{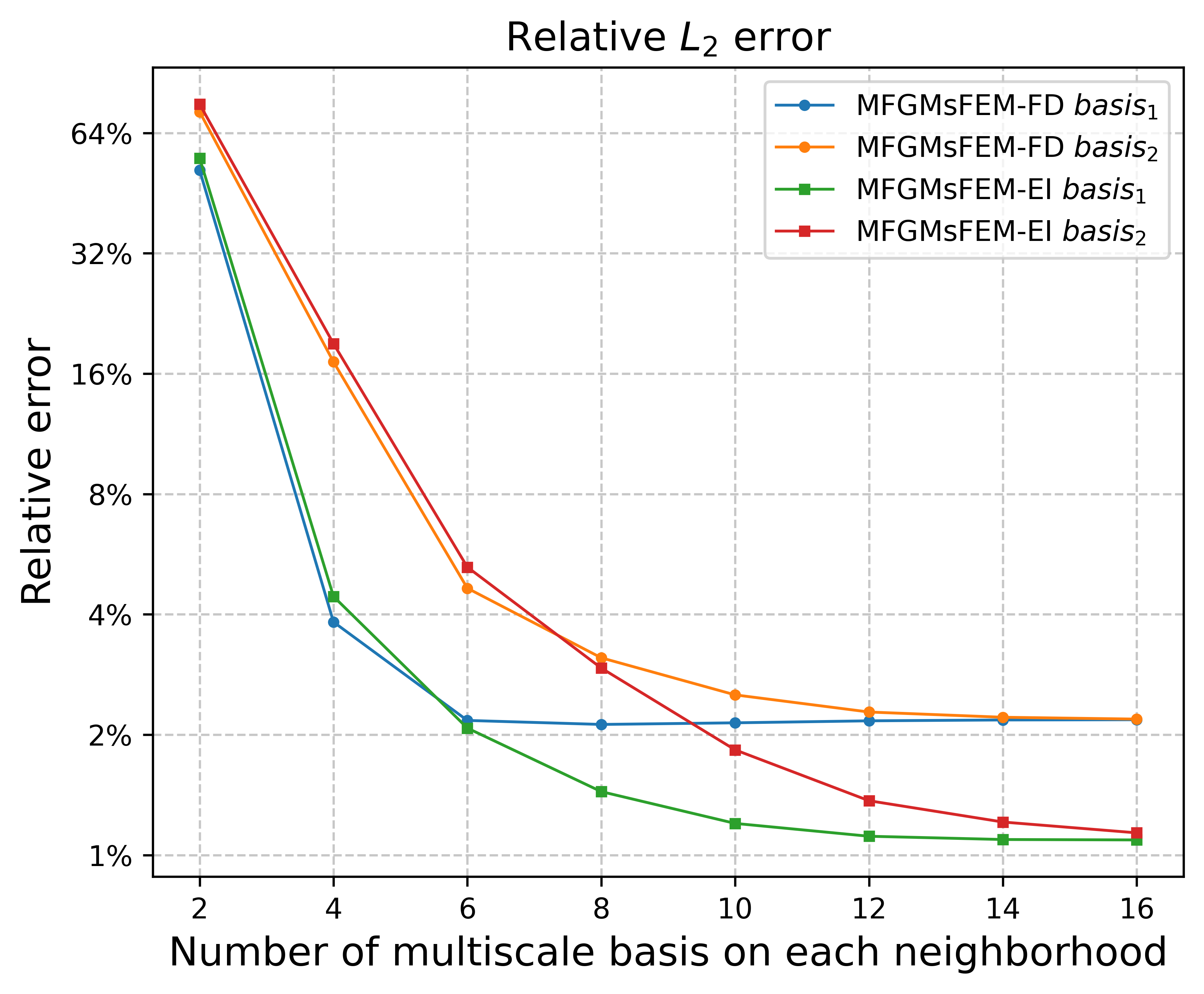}}
\end{minipage}
\begin{minipage}[h]{0.49\linewidth}
\center{\includegraphics[width=\linewidth]{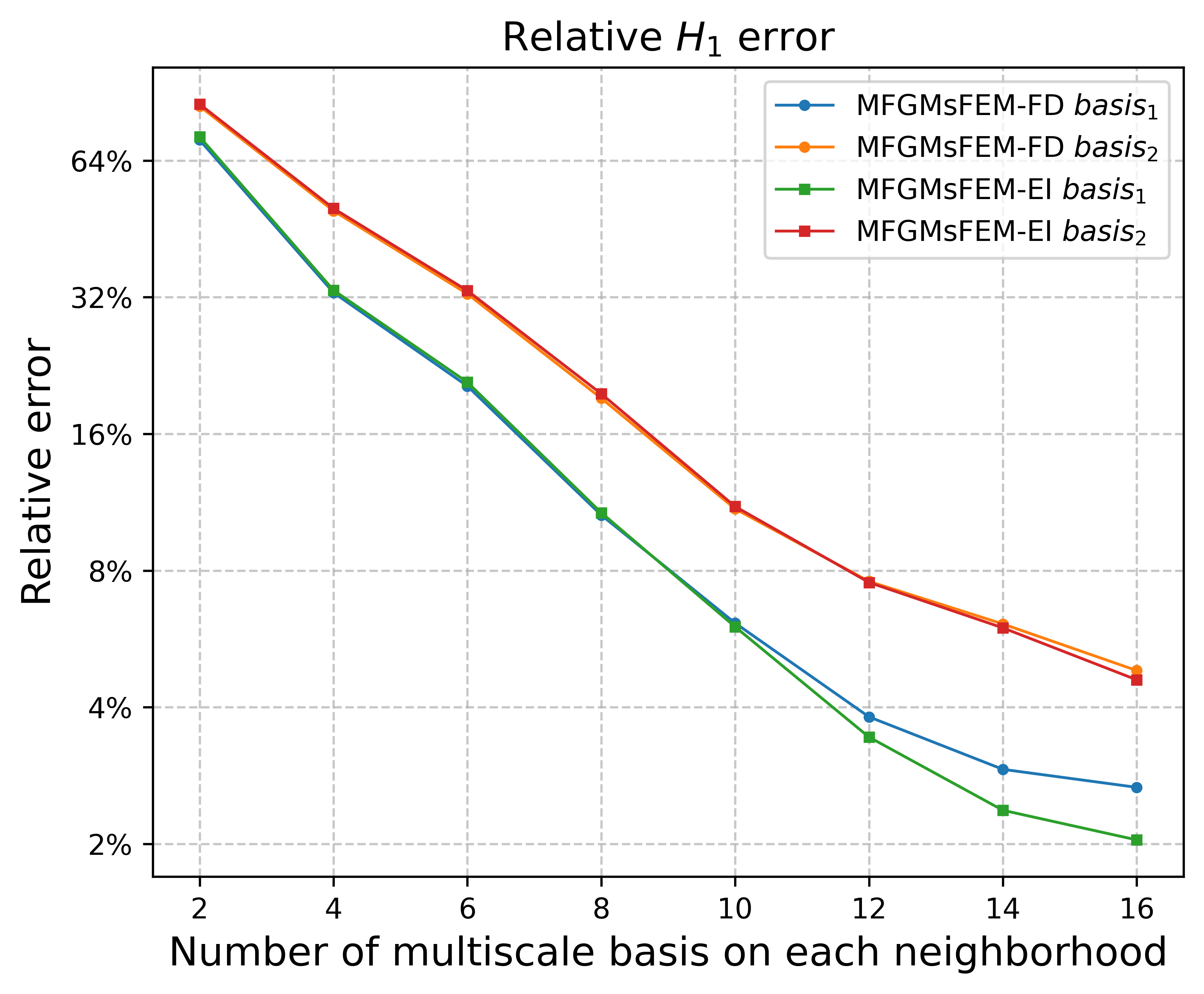}}
\end{minipage}\\
\begin{minipage}[h]{0.49\linewidth}
\center{\includegraphics[width=\linewidth]{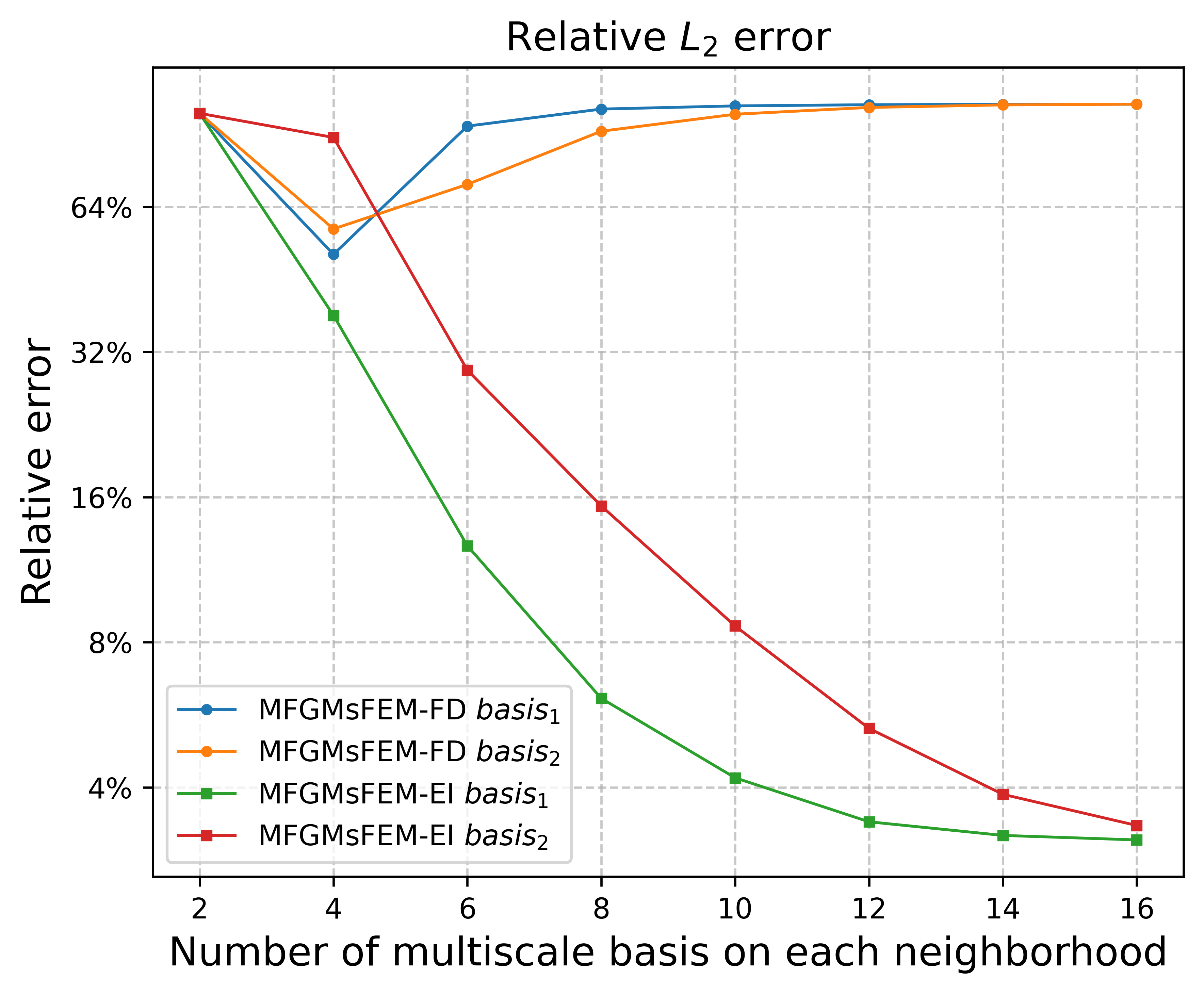}}
\end{minipage}
\begin{minipage}[h]{0.49\linewidth}
\center{\includegraphics[width=\linewidth]{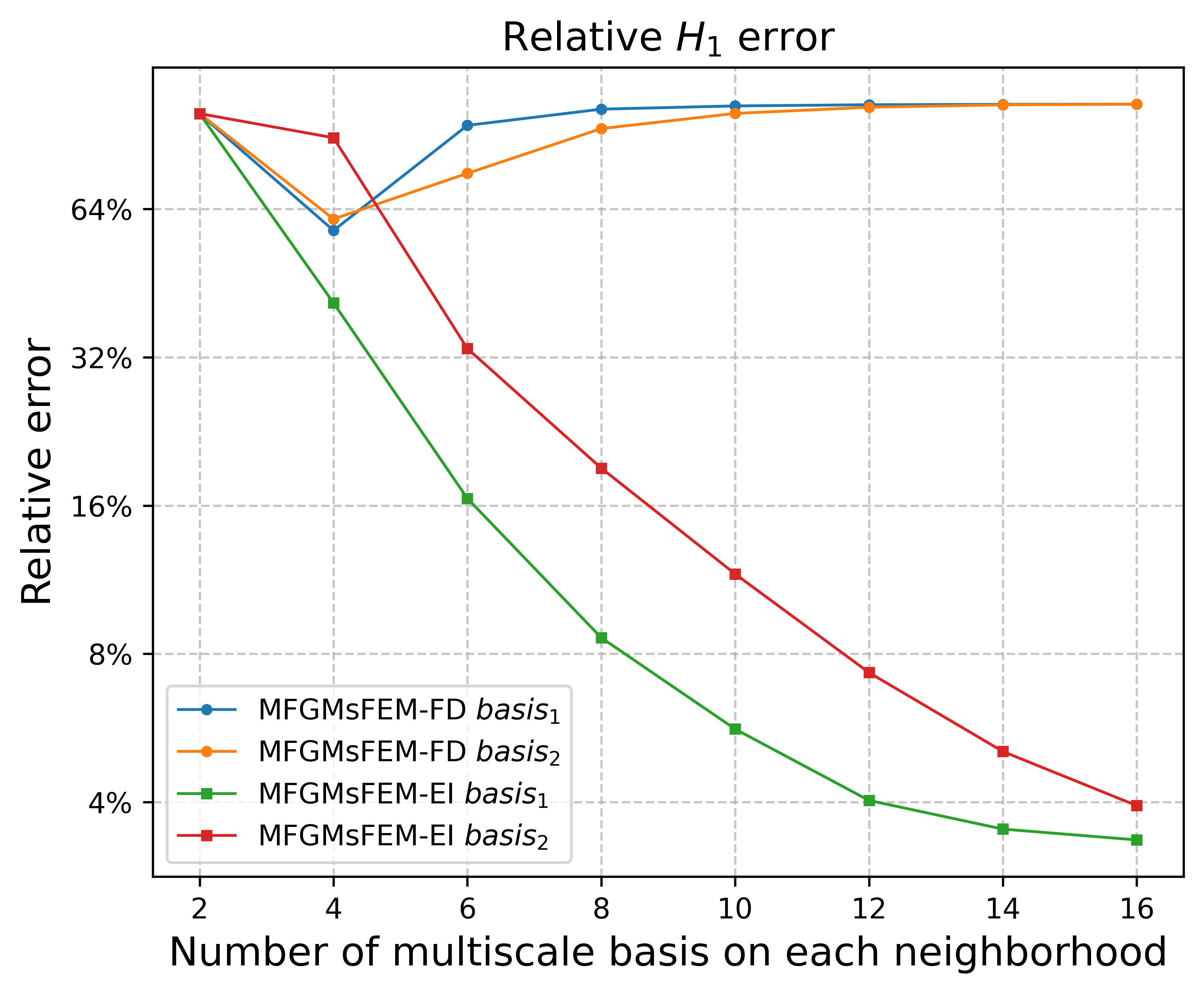}}
\end{minipage}
\end{center}
\caption{Relative errors at final time $T = 0.2$: relative $L_2$ error (left), relative $H_1$ error (right), Example 1 (top), Example 2 (bottom).}
\label{fig:error1}
\end{figure}

\begin{figure}[!htb]
\begin{center}
\begin{minipage}[h]{0.49\linewidth}
\center{\includegraphics[width=\linewidth]{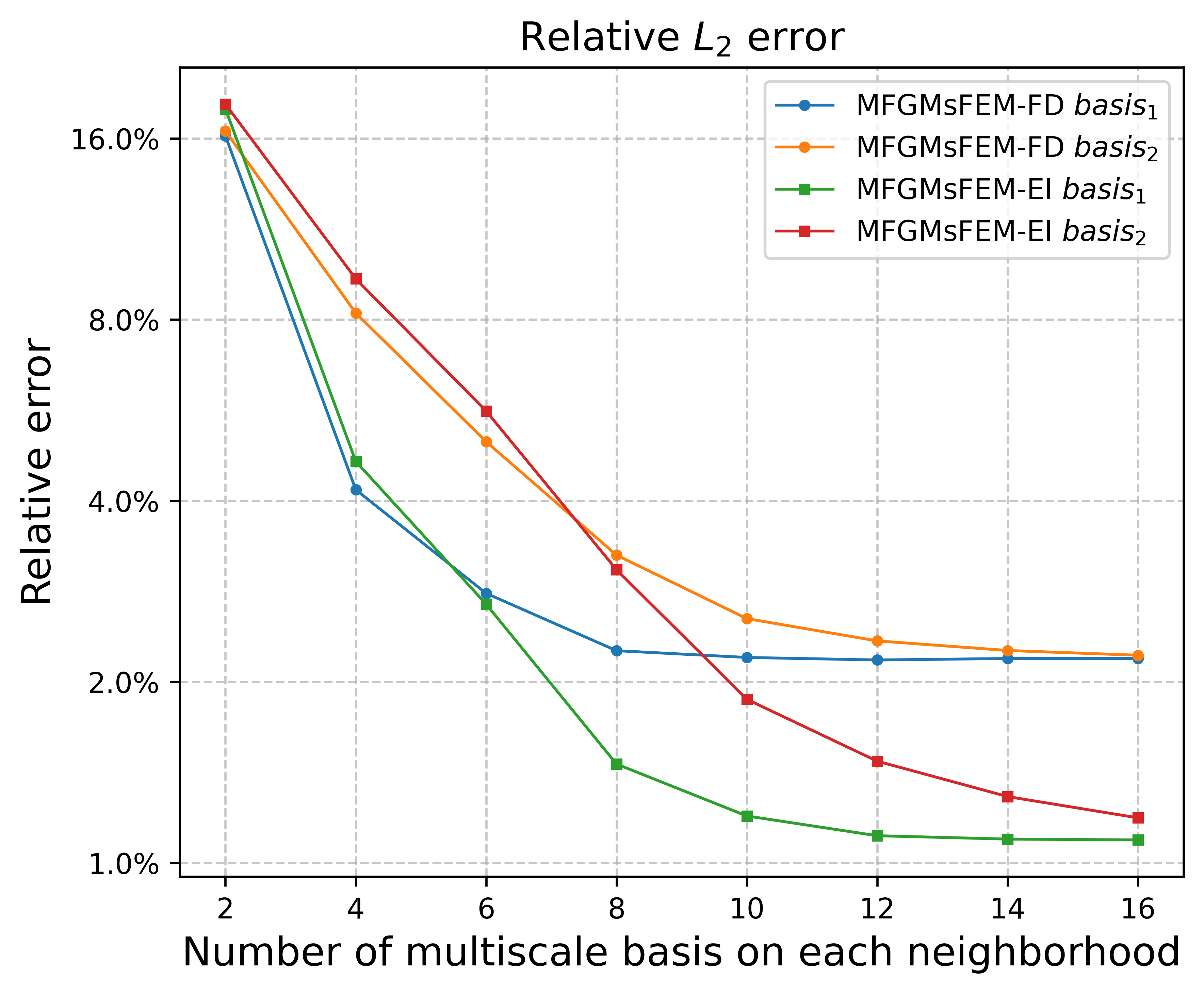}}
\end{minipage}
\begin{minipage}[h]{0.49\linewidth}
\center{\includegraphics[width=\linewidth]{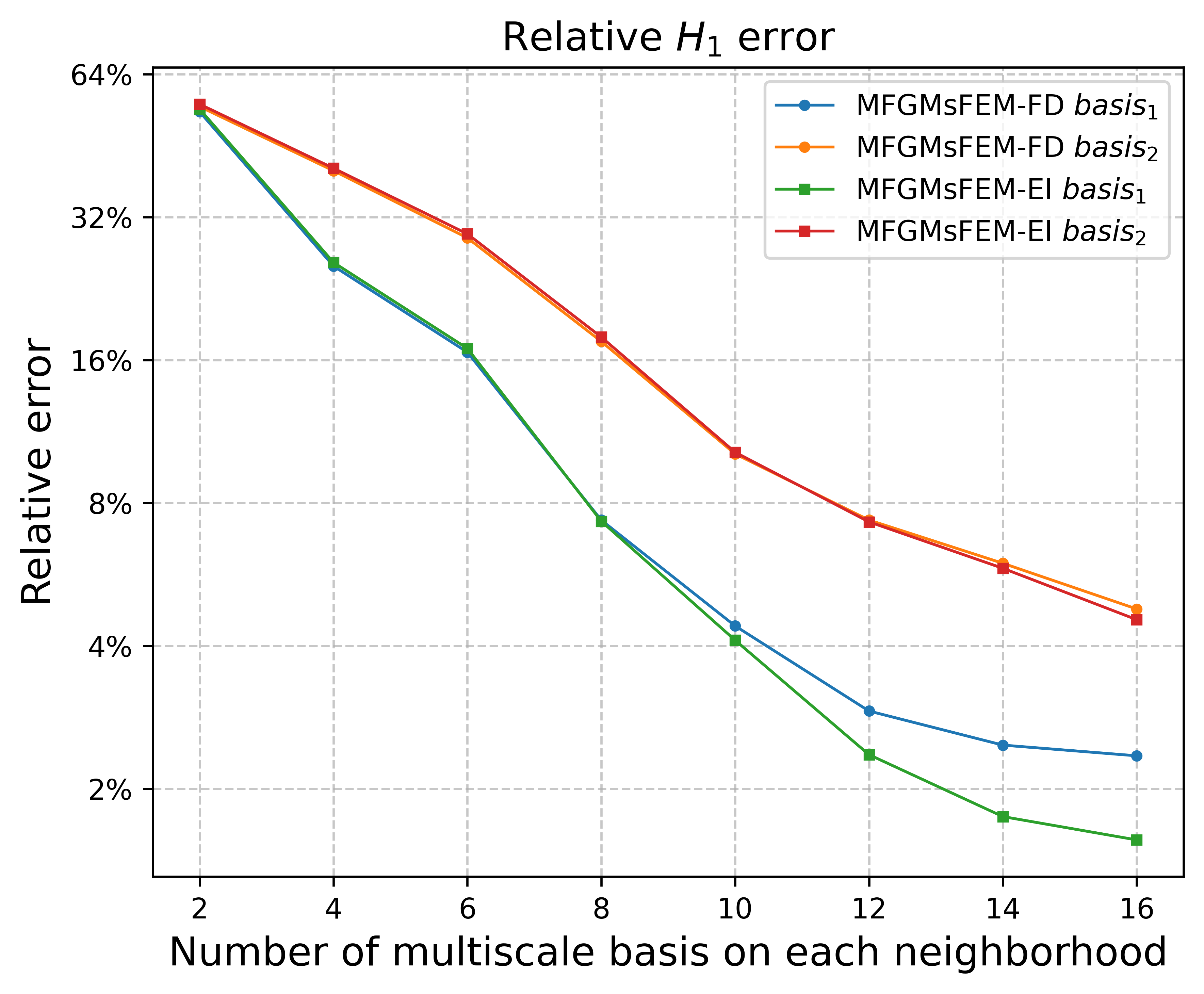}}
\end{minipage}\\
\begin{minipage}[h]{0.49\linewidth}
\center{\includegraphics[width=\linewidth]{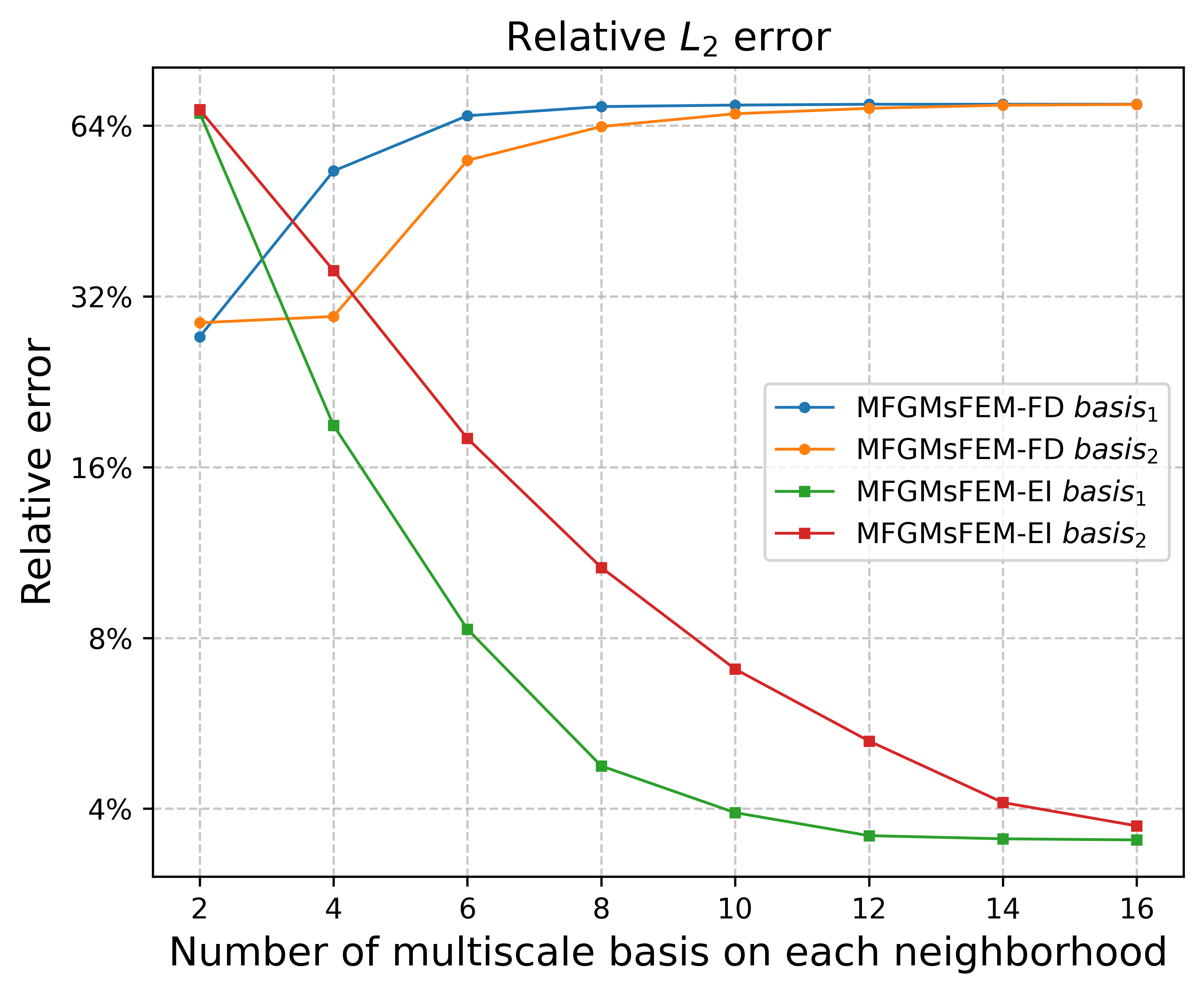}}
\end{minipage}
\begin{minipage}[h]{0.49\linewidth}
\center{\includegraphics[width=\linewidth]{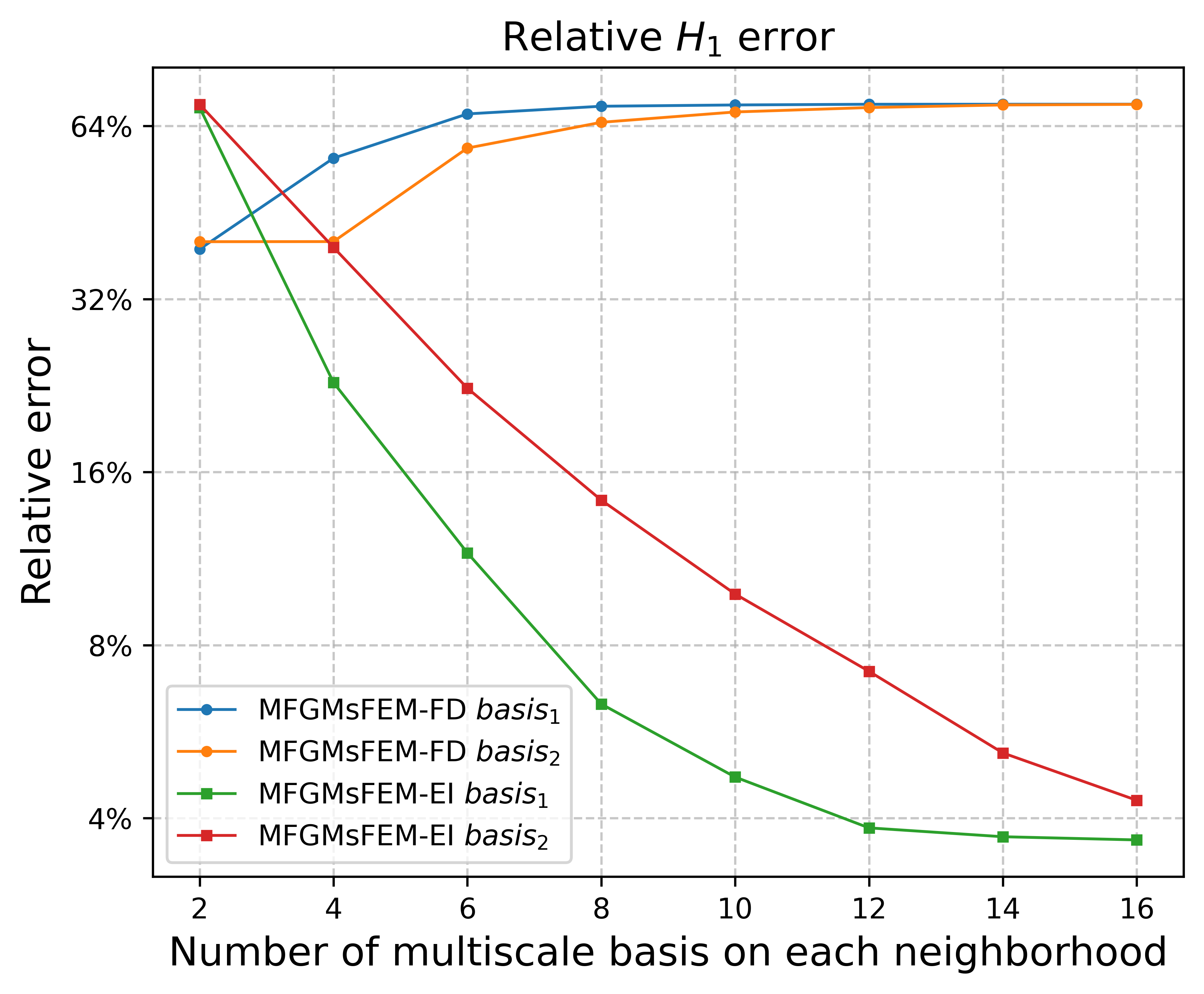}}
\end{minipage}
\end{center}
\caption{Relative errors at final time $T = 0.2$: relative $L_2$ error (left), relative $H_1$ error (right), Example 3 (top), Example 4 (bottom).}
\label{fig:error2}
\end{figure}

\begin{figure}[!htb]
\begin{center}
\begin{minipage}[h]{0.49\linewidth}
\center{\includegraphics[width=\linewidth]{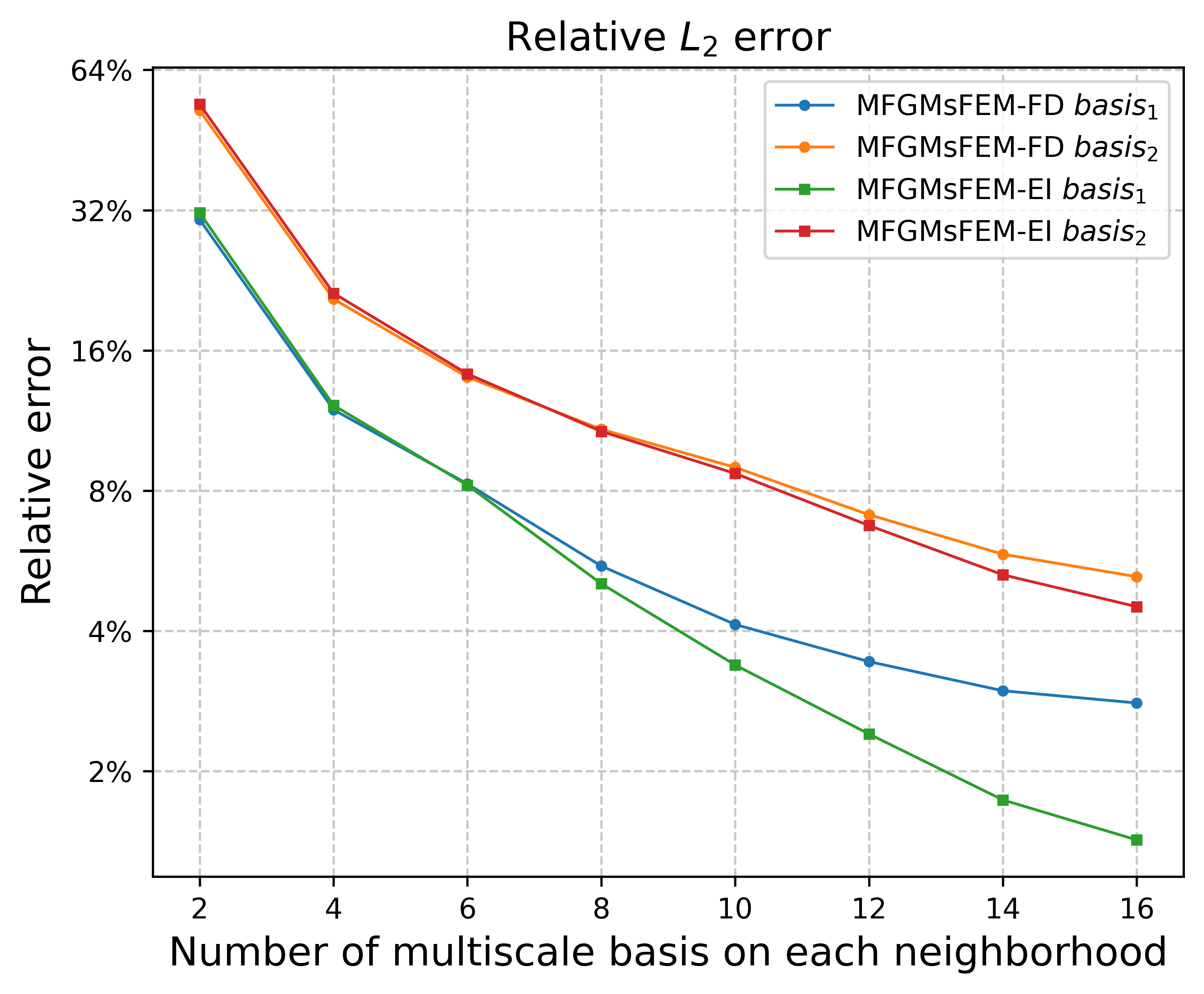}}
\end{minipage}
\begin{minipage}[h]{0.49\linewidth}
\center{\includegraphics[width=\linewidth]{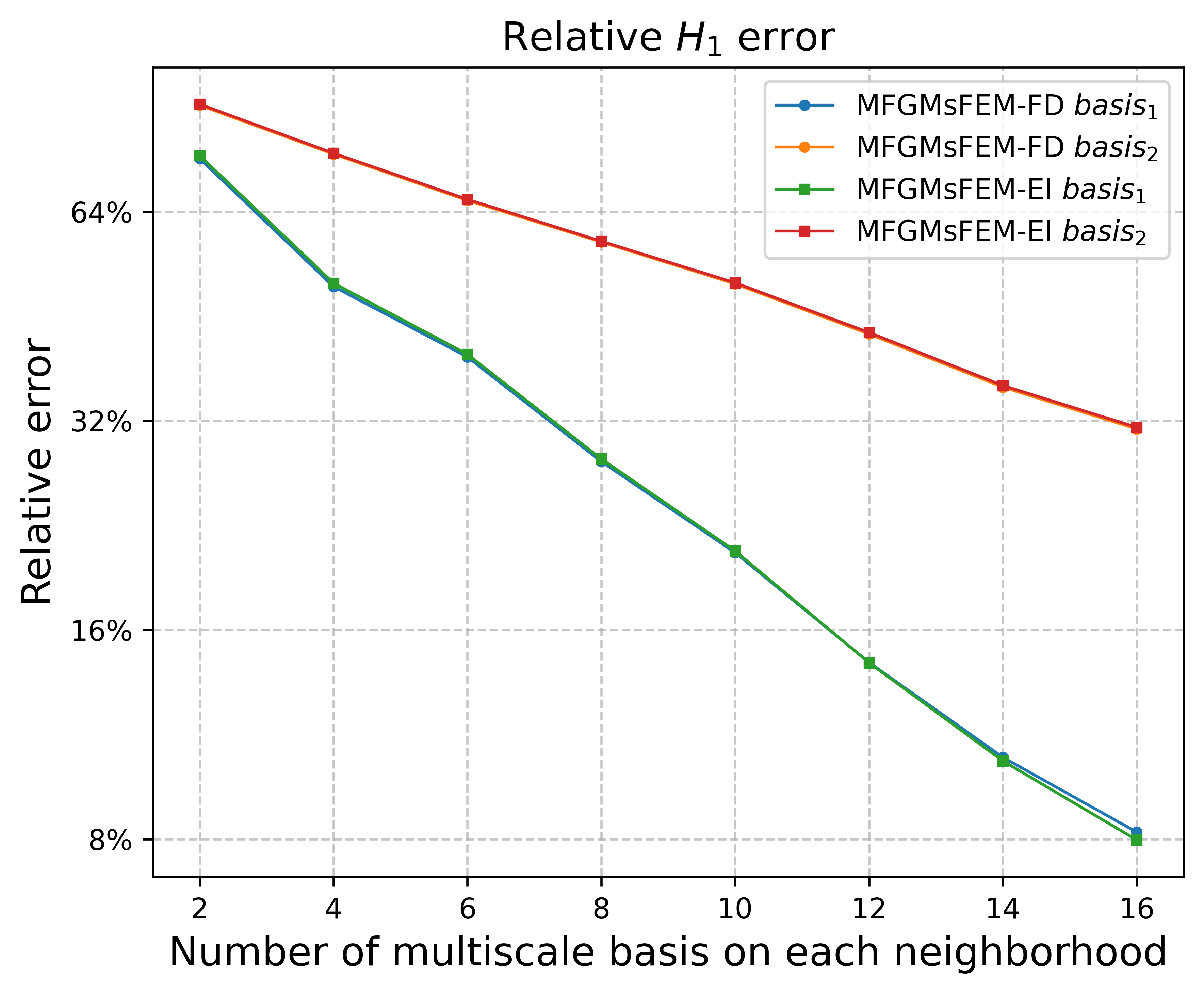}}
\end{minipage}\\
\begin{minipage}[h]{0.49\linewidth}
\center{\includegraphics[width=\linewidth]{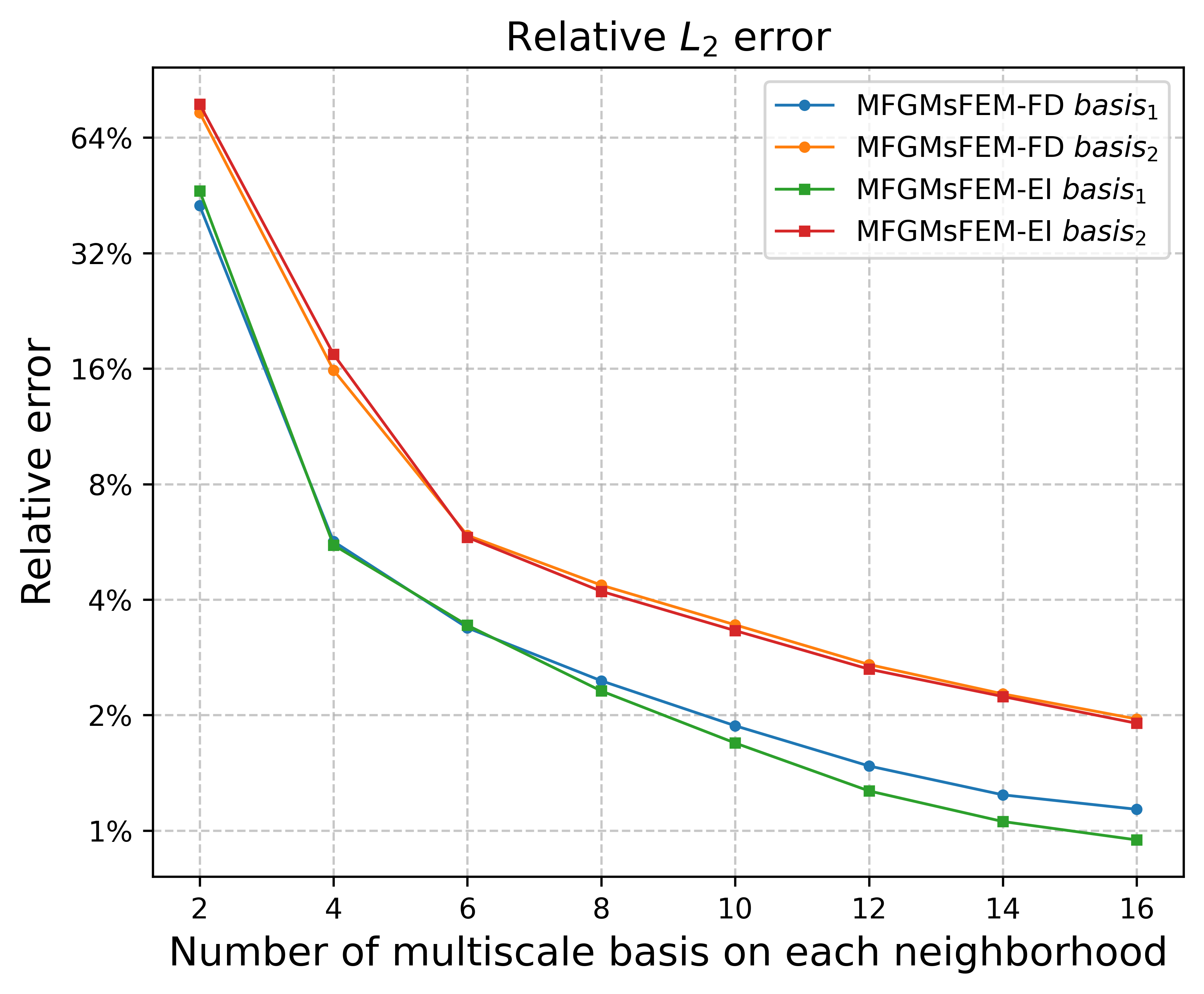}}
\end{minipage}
\begin{minipage}[h]{0.49\linewidth}
\center{\includegraphics[width=\linewidth]{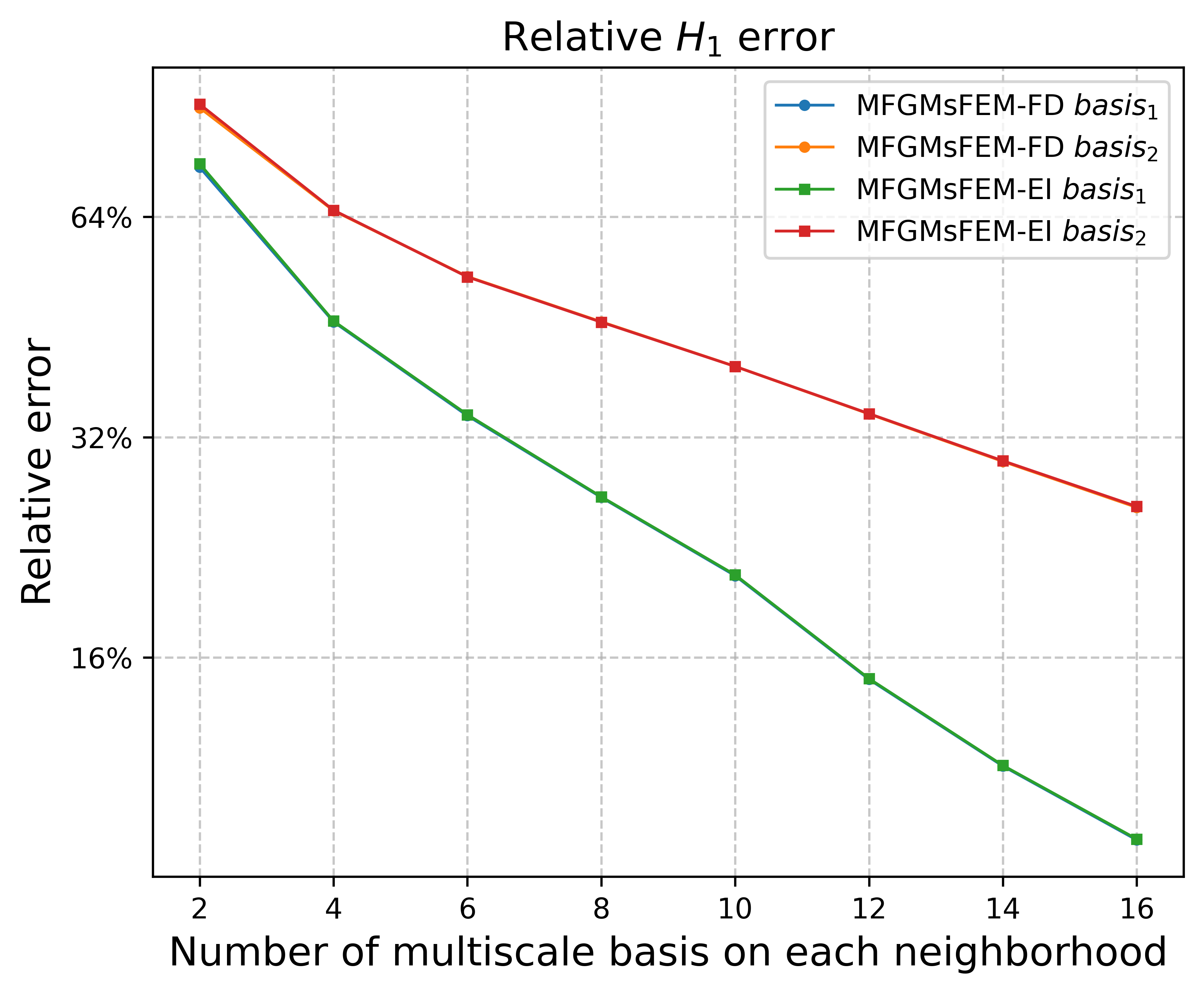}}
\end{minipage}
\end{center}
\caption{Relative errors at final time $T = 0.2$: relative $L_2$ error (left), relative $H_1$ error (right), Example 5 (top), Example 6 (bottom).}
\label{fig:error3}
\end{figure}

\begin{figure}[!htb]
\begin{center}
\begin{minipage}[h]{0.49\linewidth}
\center{\includegraphics[width=\linewidth]{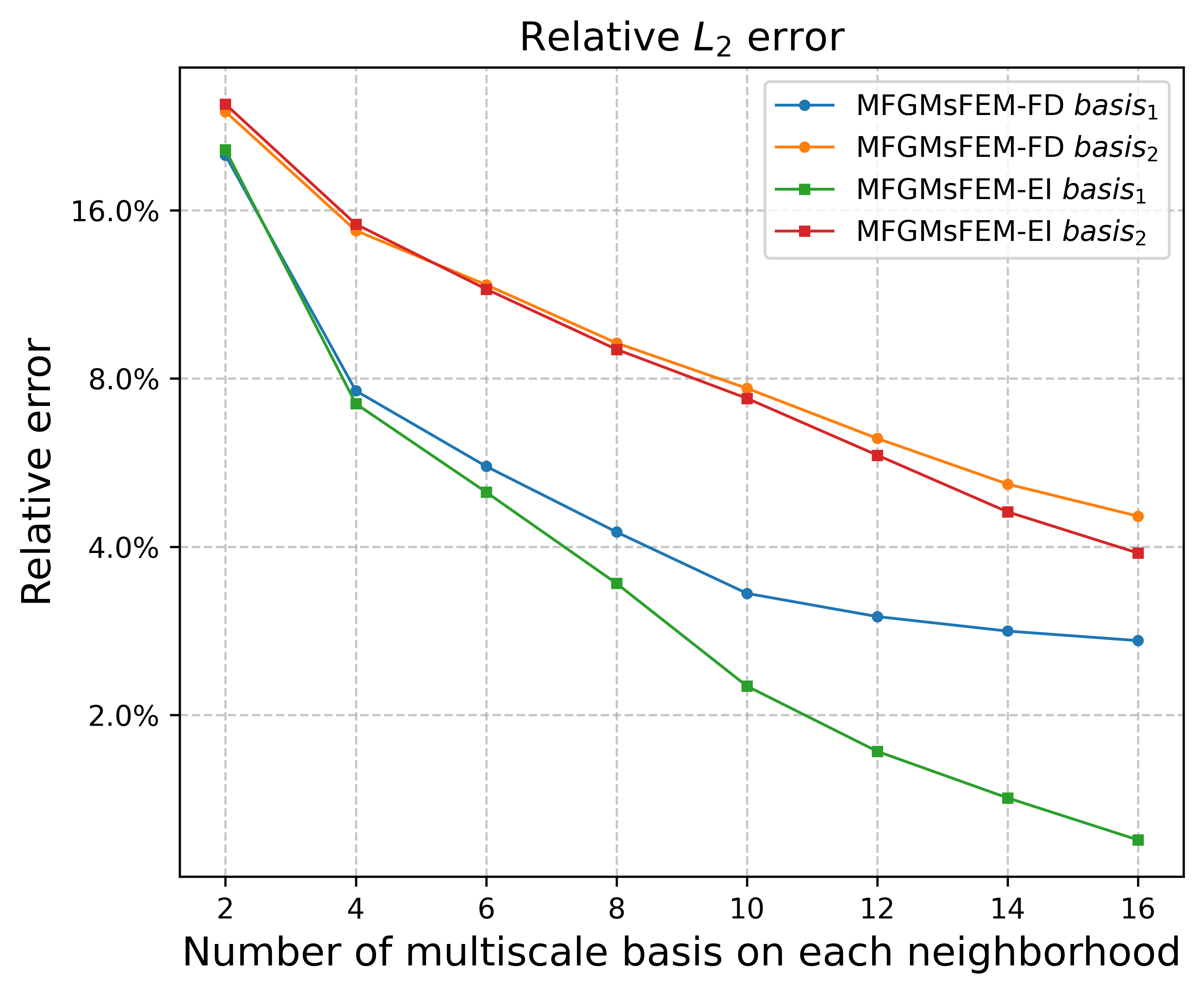}}
\end{minipage}
\begin{minipage}[h]{0.49\linewidth}
\center{\includegraphics[width=\linewidth]{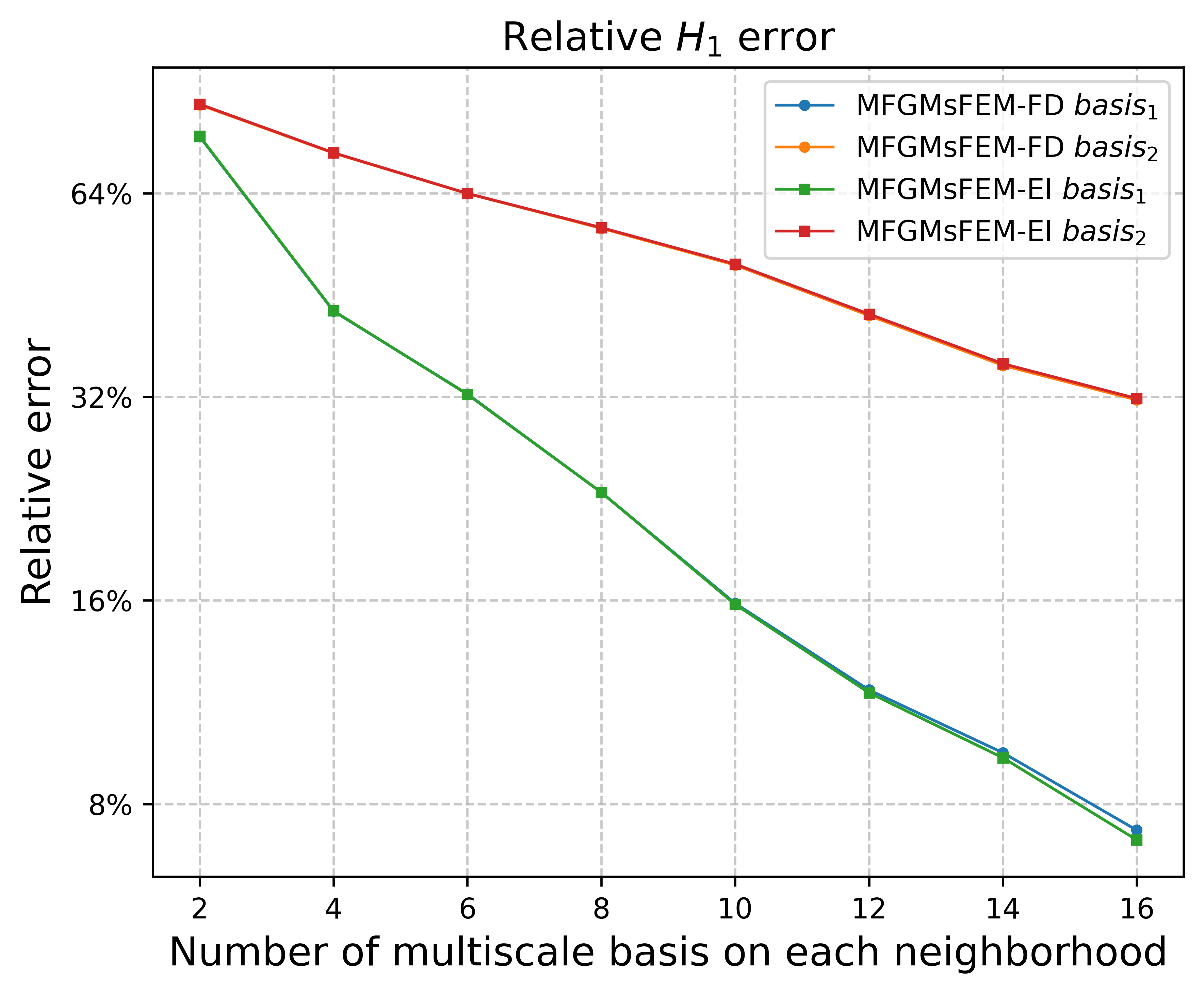}}
\end{minipage}\\
\begin{minipage}[h]{0.49\linewidth}
\center{\includegraphics[width=\linewidth]{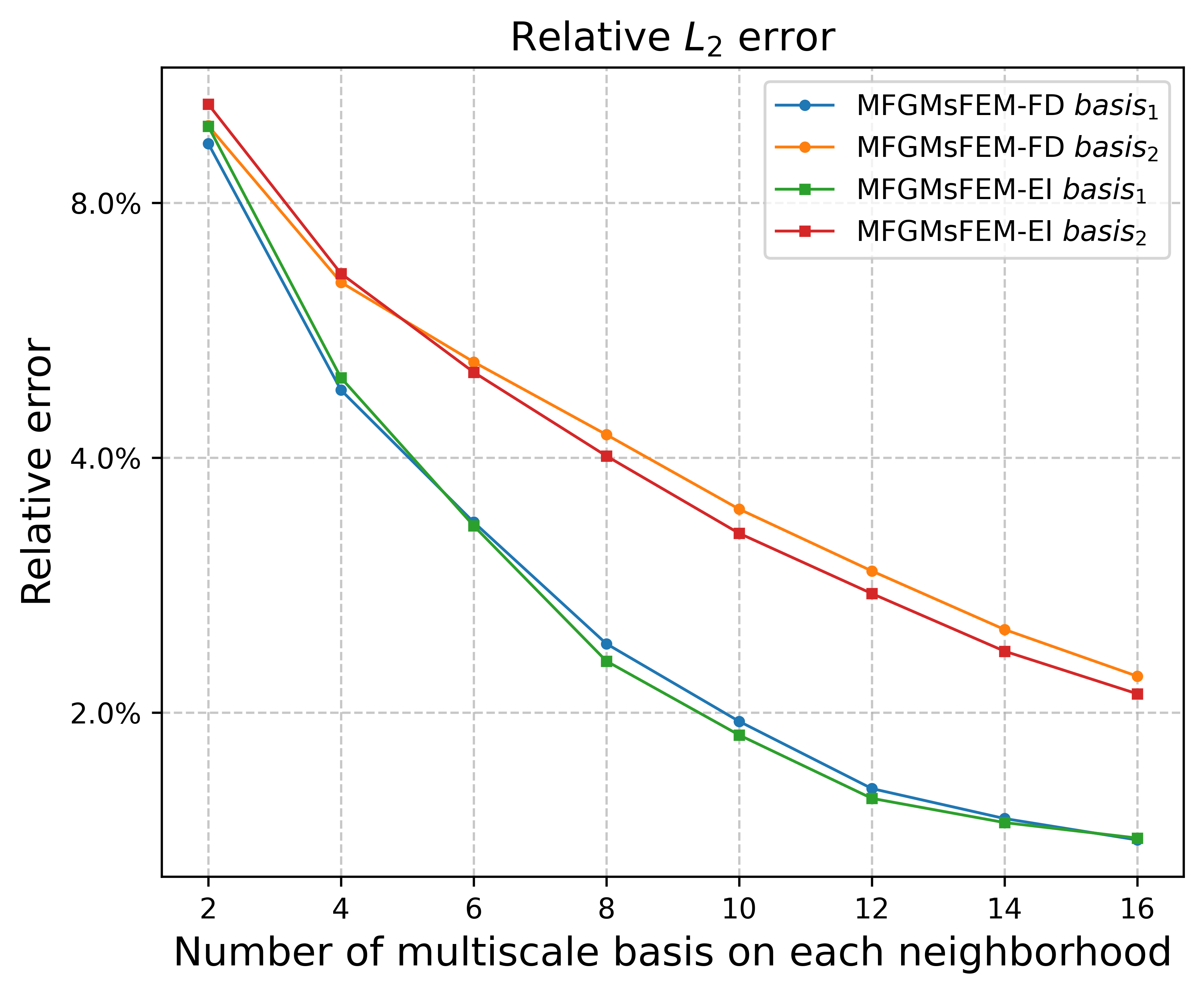}}
\end{minipage}
\begin{minipage}[h]{0.49\linewidth}
\center{\includegraphics[width=\linewidth]{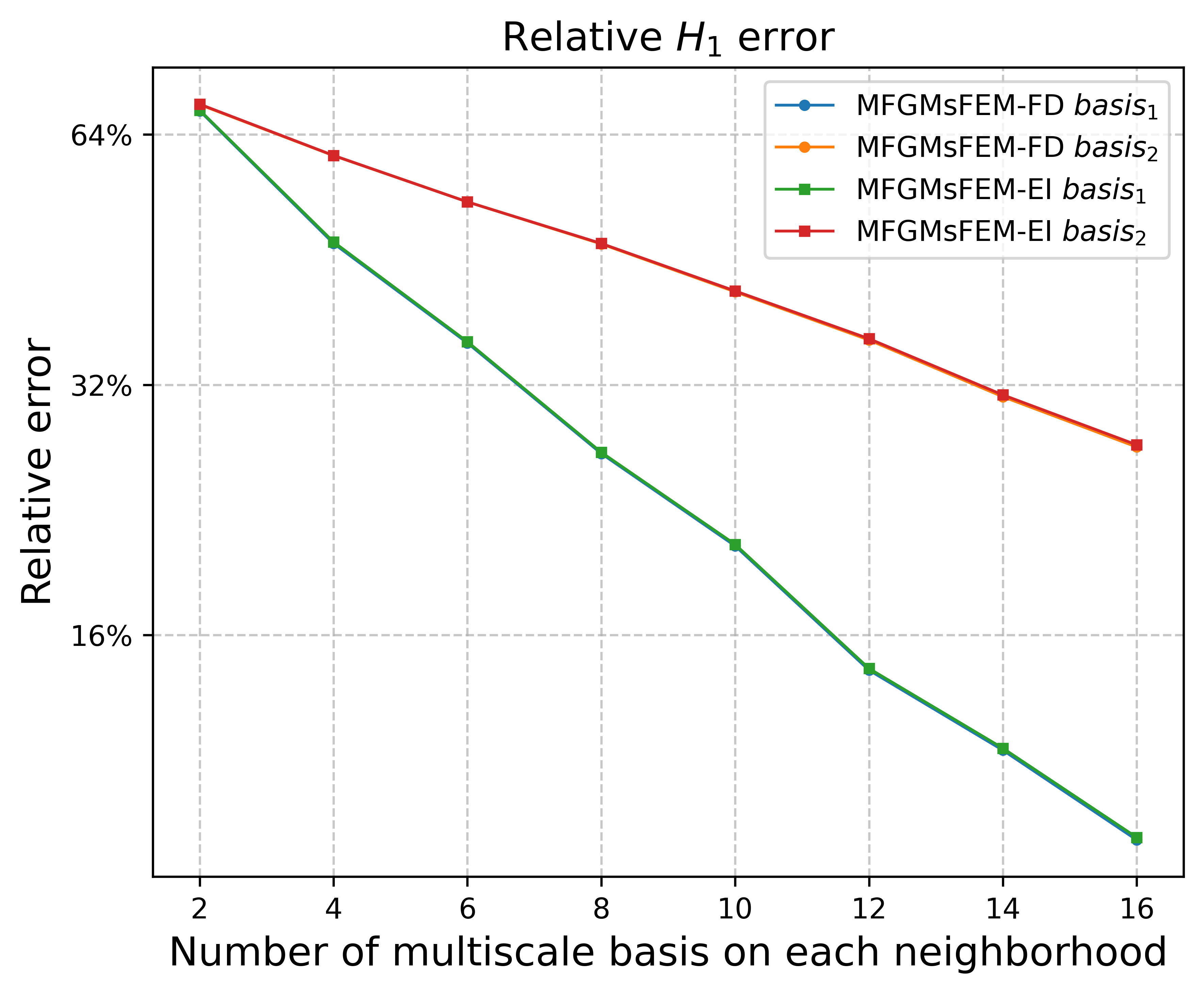}}
\end{minipage}
\end{center}
\caption{Relative errors at final time $T = 0.2$: relative $L_2$ error (left), relative $H_1$ error (right), Example 7 (top), Example 8 (bottom).}
\label{fig:error4}
\end{figure}

In Figures~4 and~5 (Examples 1--4), corresponding to $\epsilon = 1$, both $L^2$ and $H^1$ errors decrease as the number of multiscale basis functions increases, confirming the convergence of the method. The exponential integrator provides improved accuracy, especially for higher contrasts.

In Figures~6 and~7 (Examples 5--8), corresponding to $\epsilon = 1/20$, the problem becomes more challenging due to stronger advection effects. In these cases, the advantage of exponential integration is more evident, with significantly lower errors compared to the standard approach. This highlights the robustness of the proposed method in handling stiff multiscale problems.

In addition, we compare the performance of the two types of multiscale basis functions introduced in Section~3.1.2. The results indicate that both constructions provide accurate approximations; however, their behavior differs depending on the regime. The Type 1 basis functions, which incorporate advection effects only at the snapshot level, tend to provide improved accuracy in advection-dominated cases, as they better capture transport-driven features of the solution. On the other hand, the Type 2 basis functions, based on a modified spectral problem, exhibit a more uniform behavior across different regimes and remain competitive in diffusion-dominated settings. Overall, both approaches are effective, but the inclusion of advection in the basis construction becomes particularly beneficial as transport effects become more significant.
\medskip

These results demonstrate that the combination of meshfree GMsFEM and exponential integration provides an efficient and stable framework for three-dimensional multiscale advection--diffusion problems. The method preserves accuracy while allowing for larger time steps and reduced computational cost, making it a promising tool for large-scale simulations in complex heterogeneous media.

\section{Conclusions}\label{sec:conclusions}

In this work, we extended the meshfree generalized multiscale exponential integration framework introduced in \cite{nikiforov2025meshfree} to three-dimensional advection--diffusion problems. The main objective was to assess the computational viability of the method in higher-dimensional settings and in the presence of significant advection effects.

The proposed approach incorporates new constructions of multiscale basis functions tailored to advection-dominated regimes and combines them with exponential time integration to address stiffness induced by high-contrast multiscale coefficients. The numerical experiments demonstrate that the method preserves accuracy and stability while allowing for larger time steps, even in challenging three-dimensional configurations.

These results confirm that the meshfree GMsFEM combined with exponential integration is a robust and efficient methodology not only in two-dimensional settings, as previously shown in \cite{nikiforov2025meshfree}, but also for three-dimensional problems with strong transport effects. This makes the approach a promising tool for large-scale simulations in complex heterogeneous media.

\section*{Acknowledgments}
The research of Djulustan Nikiforov is supported by the Ministry of Science and Higher Education of the Russian Federation (Grant No. FSRG-2026-0009). The research of Leonardo A. Poveda is supported by Wenzhou-Kean University with an Internal Research Startup Grant (Code: ISRG2025003). The research of Dmitry Ammosov is supported by the Khalifa University Postdoctoral Fellowship Program under ADNOC's Grant No. 8434000476.

\bibliographystyle{plain}
\bibliography{references}

\begin{algorithm}[!htb]
\caption{Probabilistic construction of centroidal Voronoi points}
\label{alg:cvt}
\begin{algorithmic}[1]
\REQUIRE
    Domain $\Omega$, density $\rho(x)$, number of coarse points $N_v^H$,
    constants $\alpha_1,\alpha_2,\beta_1,\beta_2$ ($\alpha_1+\alpha_2=1$, $\beta_1+\beta_2=1$),
    samples per iteration $q$, maximum iterations $N_{\text{iter}}$.
\STATE \textbf{Initialization:}
    \STATE Generate $N_v^H$ points $\{x_i\}_{i=1}^{N_v^H}$ by random sampling according to $\rho$.
    \STATE Set counters $c_i \leftarrow 1$ for $i=1,\dots,N_v^H$.
\FOR{$k = 1$ to $N_{\text{iter}}$}
    \STATE \textbf{Random sampling:} 
        \STATE Generate $q$ points $\{y_r\}_{r=1}^{q}$ independently in $\Omega$ with probability density $\rho$.
    \STATE \textbf{Nearest‑neighbor assignment:}
        \FOR{$r = 1$ to $q$}
            \STATE $i(r) \leftarrow \arg\min_{i} \|y_r - x_i\|$ \COMMENT{closest coarse point}
            \STATE Add $y_r$ to the collection $C_{i(r)}$.
        \ENDFOR
    \STATE \textbf{Update coarse points:}
        \FOR{$i = 1$ to $N_v^H$}
            \IF{$C_i \neq \emptyset$}
                \STATE $\bar{x}_i \leftarrow \frac{1}{|C_i|}\sum_{y \in C_i} y$ \COMMENT{centroid of assigned samples}
                \STATE $x_i \leftarrow \frac{(\alpha_1 c_i + \beta_1)\, x_i + (\alpha_2 c_i + \beta_2)\, \bar{x}_i}{c_i + 1}$
                \STATE $c_i \leftarrow c_i + 1$
            \ENDIF
            \STATE Clear $C_i$ for the next iteration.
        \ENDFOR
\ENDFOR
\ENSURE The final coarse point set $\{x_i\}_{i=1}^{N_v^H}$.
\end{algorithmic}
\end{algorithm}

\begin{algorithm}[!htb]
\caption{Determination of coarse-element radii}
\label{alg:radii}
\begin{algorithmic}[1]
\REQUIRE
    Domain $\Omega$, coarse points $\{x_i\}_{i=1}^{N_v^H}$,
    a set of pseudo‑points $\mathcal{P}$ that densely samples $\Omega$ (e.g., a uniform grid),
    overlap factor $\gamma > 1$.
\STATE \textbf{Initialization:} Set $r_i \leftarrow 0$ for all $i$.
\FOR{each $p \in \mathcal{P}$}
    \STATE $i(p) \leftarrow \arg\min_{i} \|p - x_i\|$ \COMMENT{closest coarse point}
    \STATE $r_{i(p)} \leftarrow \max\bigl(r_{i(p)},\; \|p - x_{i(p)}\|\bigr)$
\ENDFOR
\FOR{$i = 1$ to $N_v^H$}
    \STATE $r_i \leftarrow \gamma \cdot r_i$ \COMMENT{enlarge to ensure overlap}
\ENDFOR
\ENSURE Radii $\{r_i\}_{i=1}^{N_v^H}$ such that $\Omega \subset \bigcup_i S_i$ with $S_i = B(x_i, r_i)$.
\end{algorithmic}
\end{algorithm}
\end{document}